\documentclass[10pt, reqno]{amsart}

\usepackage[utf8]{inputenc}
\usepackage{amsmath}
\usepackage{amssymb}
\usepackage{mathrsfs}
\usepackage{graphics}
\usepackage[dvipdf]{graphicx}
\usepackage{geometry}
\usepackage{setspace}
\usepackage{color}
\usepackage{amsthm}
\usepackage{paralist}
\usepackage{dsfont}
\usepackage{verbatim}
\usepackage[numbers]{natbib}
\usepackage[breaklinks, colorlinks, linkcolor = black, citecolor = black, filecolor = black, urlcolor = black, bookmarksopen=true]{hyperref} 

\theoremstyle{definition}
\newtheorem{definition}{Definition}[section]

\newtheorem{example}[definition]{Example}
\theoremstyle{plain}
\newtheorem{lemma}[definition]{Lemma}
\newtheorem{thm}[definition]{Theorem}
\newtheorem{prop}[definition]{Proposition}

\numberwithin{equation}{section}

\allowdisplaybreaks

\newcommand{\auxil}{\mathcal{A}}
\newcommand{\metr}{\mathbf{d}}

\newcommand{\aV}{\mathcal{V}}
\newcommand{\ent}{\mathcal{E}}

\newcommand{\dd}{\,\mathrm{d}}

\newcommand{\dv}{\mathrm{div}\,}
\newcommand{\dff}{\mathrm{D}}
\newcommand{\prb}{\mathscr{P}}

\newcommand{\ball}{\mathbb{B}}
\newcommand{\calH}{{\mathcal{H}}}

\newcommand{\eps}{\varepsilon}
\newcommand{\W}{\mathbf{W}}

\newcommand{\X}{\mathbf{X}}

\newcommand{\N}{{\mathbb{N}}}
\newcommand{\R}{{\mathbb{R}}}
\newcommand{\Rd}{{\R^d}}
\newcommand{\flow}{\mathsf{S}}
\newcommand{\eins}[1]{\mathbf{1}_{#1}}
\newcommand{\Id}{\mathds{1}}

\newcommand{\tT }{\mathrm{T}}

\newcommand{\trc}{\mathrm{tr}}

\newcommand{\supp}{\mathrm{supp}\,}

\newcommand{\einsvec}{\mathsf{e}}

\newcommand{\mom}[1]{\mathfrak{m}_2(#1)}

\newcommand{\gau}[1]{\left\lfloor#1\right\rfloor}
\newcommand{\ban}{\mathbf{Y}}
\newcommand{\dom}{\mathrm{Dom}}
\newcommand{\argmin}{\operatornamewithlimits{argmin}}

\renewcommand{\tilde}{\widetilde}
\renewcommand{\bar}{\overline}

\begin{document}

\begin{abstract}
We prove the existence of nonnegative weak solutions to a class of second and fourth order nonautonomous nonlinear evolution equations with an explicitly time-dependent mobility function posed on the whole space $\Rd$, for arbitrary $d\ge 1$. Exploiting a very formal gradient flow structure, the cornerstone of our proof is a modified version of the classical minimizing movement scheme for gradient flows. The mobility function is required to satisfy---at each time point separately---the conditions by which one can define a modified Wasserstein distance on the space of probability densities with finite second moment. The explicit dependency on the time variable is assumed to be at least of Lipschitz regularity. We also sketch possible extensions of our result to the case of bounded spatial domains and more general mobility functions.
\end{abstract}

\title[Time-dependent mobility]{Well-posedness of evolution equations with time-dependent nonlinear mobility: a modified minimizing movement scheme}
\author[Jonathan Zinsl]{Jonathan Zinsl}
\address{Zentrum f\"ur Mathematik \\ Technische Universit\"at M\"unchen \\ 85747 Garching, Germany}
\email{zinsl@ma.tum.de}
\keywords{Non-autonomous equation; gradient flow; nonlinear mobility; modified Wasserstein distance; minimizing movement scheme}
\date{\today}
\subjclass[2010]{Primary: 35K30; Secondary: 35A15, 35D30.}
\thanks{This research has been supported by the DFG Collaborative Research Center TRR 109 ``Discretization in Geometry and Dynamics''}

\maketitle

\section{Introduction}\label{sec:intro}
In this article, the following nonautonomous partial differential equation is considered:
\begin{align}\label{eq:pde}
\partial_t u(t,x)&=\dv\left(m(t,u(t,x))\nabla\frac{\delta \ent}{\delta u}(u(t,x))\right),\qquad t>0,\,x\in\Rd,
\end{align}
together with the initial condition $u(0,\cdot)=u_0$ on $\Rd$ for $d\ge 1$. In addition, we seek for a \emph{nonnegative} solution $u(t,x)\ge 0$ for all $t\ge 0$ and almost all $x\in\Rd$.

Above, $\frac{\delta \ent}{\delta u}$ denotes the first variation of the \emph{free energy} functional $\ent$ in $L^2$. Our assumptions on $\ent$ are specified below in the conditions (E1) and (E2), respectively. There, we distinguish two main classes of energy functionals. First, in the case (E1), $\ent$ is assumed to be of the form
\begin{align}\label{eq:ent2o}
\ent(u)=\int_\Rd \left[f(u(x))+\phi(x)u(x)\right]\dd x,
\end{align}
where the internal energy density $f$ is uniformly convex in $u$, and the \emph{confinement potential} $\phi$ is bounded from below and grows at most quadratically as $|x|\to\infty$. Equation \eqref{eq:pde} then reads as a nonlinear second order drift-diffusion equation,
\begin{align}\label{eq:pde2o}
\partial_t u&=\dv\left[m(t,u)(f''(u)\nabla u+\nabla \phi)\right].
\end{align}
One may also allow for \emph{gradient-dependent} energy density: 
\begin{align}\label{eq:ent4o}
\ent(u)=\int_\Rd \left[f(\nabla u(x),u(x))+\phi(x)u(x)\right]\dd x,
\end{align}
where again $f$ is assumed to satisfy a certain convexity condition. Equation \eqref{eq:pde} then is the fourth-order equation
\begin{align}\label{eq:pde4o}
\partial_t u&=\dv\left[m(t,u)\nabla(-\dv\,\nabla_pf(\nabla u,u)+\partial_zf(\nabla u,u)+\phi)\right].
\end{align}
We do not seek for full generality of energy functionals $\ent$ here, since our main interest is focussed on the \emph{time-dependent mobility function} $m:\R_{\ge 0}\times \R_{\ge 0}\to \R_{\ge 0}$ which turns \eqref{eq:pde} into a \emph{non-autonomous} nonlinear evolution equation.

When $m$ does not explicitly depend on $t$, it is known that \eqref{eq:pde} possesses a variational structure \cite{dns2009,lisini2010,lisini2012}, if $m$ is nonnegative and concave on the interior of an interval $[0,S]$, the so-called \emph{value space} \cite{zm2014}, with $S\in\R_{\ge 0}\cup\{+\infty\}$. In this work, we require that at each fixed time $t\ge 0$, $m(t,\cdot)$ is an admissible mobility in the sense of \cite{dns2009,lisini2010} (see condition (M1) below), where the corresponding value spaces $[0,S(t)]$ are assumed to be expanding over time, i.e., $S$ is nondecreasing. For the dependency on the time variable, certain regularity---at least Lipschitz---conditions are needed. The detailed assumptions on $m$ are presented in (M1)--(M4) below. 
\begin{example}[Paradigmatic examples]
The following mobility functions admit the conditions (M1)--(M4) below:
\begin{enumerate}[(a)]
\item \emph{Finite value spaces}: For all $t\ge 0$, put
\begin{align*}
m(t,z)=z(S(t)-z)\qquad\text{for each }z\in [0,S(t)],
\end{align*}
for some sufficiently regular nondecreasing function $S$. Choosing the quadratic free energy
\begin{align}\label{eq:entQ}
\ent_Q(u)=\frac12\int_\Rd |\nabla u(x)|^2\dd x+\frac12 \int_\Rd u(x)^2\dd x,
\end{align}
which is admissible in (E2), the fourth-order equation \eqref{eq:pde4o} reads as a variant of the \emph{Cahn-Hilliard} equation,
\begin{align*}
\partial_t u&=-\dv(u(S(t)-u)\nabla\Delta u)+\dv(u(S(t)-u)\Delta u).
\end{align*}
\item \emph{Infinite value space}, i.e., $S\equiv +\infty$: Let $\eps>0$ and let $\alpha:\R_{\ge 0}\to(0,1]$ be sufficiently regular. Set
\begin{align*}
m(t,z)=(z+\eps)^{\alpha(t)}-\eps^{\alpha(t)}\quad\text{for all }z\ge 0.
\end{align*}
With the quadratic free energy $\ent_Q$ from \eqref{eq:entQ}, \eqref{eq:pde4o} is a perturbed version of the \emph{thin-film} equation,
\begin{align*}
\partial_t u&=-\dv\big(u((u+\eps)^{\alpha(t)})-\eps^{\alpha(t)})\nabla\Delta u\big)+\dv\big(u((u+\eps)^{\alpha(t)}-\eps^{\alpha(t)})\Delta u\big).
\end{align*}
\end{enumerate}
Notice that in (b), the case $\eps=0$ is (whenever $\alpha$ is not identically equal to $1$) not allowed since condition (M3) requires Lipschitz-continuity of $m(t,\cdot)$ also in the state variable $z$. For a possible generalization where (M3) can be dispensed of, see Section \ref{sec:gen}.
\end{example}
We use the (pointwise in $t$) formal gradient structure of \eqref{eq:pde} from \cite{dns2009,lisini2010} on the space
\begin{align}\label{eq:spacet}
\X(t):=\left\{u\in L^1(\Rd):~0\le u(x)\le S(t)\text{ a.e. on }\Rd,~\|u\|_{L^1}=1\text{ and }\int_\Rd |x|^2 u(x)\dd x<\infty\right\},
\end{align}
with respect to the metric $\W_{m(t,\cdot)}$ induced by $m(t,\cdot)$:
\begin{align}\label{eq:metr}
\begin{split}
&\W_{m(t,\cdot)}(u_0,u_1)\\
&\quad=\inf\left\{\int_0^1\int_\Rd \frac{|w_s|^2}{m(t,u_s)}\dd x \dd s:~(u_s,w_s)_{s\in[0,1]}\in \mathscr{C},~u_s|_{s=0}=u_0,~u_s|_{s=1}=u_1\right\}^{1/2},
\end{split}
\end{align}
where $\mathscr{C}$ is a suitable subclass of solutions to the \emph{continuity equation} $\partial_s u_s=-\dv w_s$ in the sense of distributions on $[0,1]\times\Rd$ (see \cite{dns2009,lisini2010} for more details). By convention, we set $\W_{m(t,\cdot)}(u_0,u_1)=+\infty$ if $u_0$ and $u_1$ are not elements of $\X(t)$ both. If $m(t,\cdot)$ is a linear function, we recover (a scalar multiple of) the classical $L^2$-Wasserstein distance $\W_2$ on the space of probability measures \cite{brenier2000} (see Section \ref{sec:pre} below).

Specifically, we prove the existence of nonnegative weak solutions to \eqref{eq:pde} using a modified version of the classical \emph{minimizing movement scheme} for gradient flows (which has been employed for various evolution equations with metric gradient flow structure, e.g. \cite{jko1998,otto2001,savare2008,gianazza2009,matthes2009,laurencot2011,lisini2012,zinsl2012,zm2014,blanchet2014}, also with spatially varying coefficients \cite{petrelli2004,lisini2009,loibl2015}):
\begin{definition}[Modified minimizing movement scheme]\label{def:minmov}
Define the space 
\begin{align}\label{eq:space}
\X:=\left\{u\in L^1(\Rd):~0\le u(x)\text{ a.e. on }\Rd,~\|u\|_{L^1}=1\text{ and }\mom{u}<\infty\right\},
\end{align}
where $\mom{u}=\int_\Rd |x|^2 u(x)\dd x$ denotes the second moment of $u$.\\
Let now $\tau>0$ and $u_0\in\X$ be given. A sequence $(u_\tau^n)_{n\in\N}$ is obtained via the \emph{modified minimizing movement scheme} if $u_\tau^0=u_0$ and for all $n\in\N$,
\begin{align}\label{eq:mms}
u_\tau^n\in \argmin_{u\in\X}\left[\frac1{2\tau}\W_{m(n\tau,\cdot)}(u_\tau^{n-1},u)^2+\ent(u)\right].
\end{align}
The \emph{discrete solution} $u_\tau:\R_{\ge 0}\times\X\to \R_{\ge 0}\cup\{+\infty\}$ is defined by piecewise constant interpolation along the sequence $(u_\tau^n)_{n\in\N}$, that is
\begin{align*}
u_\tau(t):=u_\tau^n,\quad\text{for }n:=\left\lceil\frac{t}{\tau}\right\rceil\quad\text{for all }t\ge 0.
\end{align*}
\end{definition}
We now present our main assumptions and results.
\subsection{Assumptions and main results}\label{subsec:main}
For the mobility function, we require the following.
\begin{definition}(Admissible mobility functions)\label{subsec:mobility}~
\begin{enumerate}[({M}1)]
\item There exists a nondecreasing function $S:\R_{\ge 0}\to \R_{>0}\cup\{+\infty\}$ such that for fixed $t\ge 0$, the map $m(t,\cdot)\in C^2((0,S(t)))$ is an admissible mobility for the definition of a distance $W_{m(t,\cdot)}$, viz.
\begin{align*}
m(t,z)&>0\quad\text{for all } z\in(0,S(t)),\quad\text{and}\quad m(t,z)=0\quad\text{for all }z\in \R_{\ge 0}\setminus (0,S(t)),\\
\partial_{z}^2 m(t,z)&\le 0\quad\text{for all } z\in(0,S(t)).
\end{align*}
\item For fixed $z\ge 0$, the map $m(\cdot,z)$ is locally Lipschitz-continuous on $\R_{\ge 0}$.
\item There exist $M_1,M_2\in C(\R_{\ge 0};\R_{>0})$ such that for all $t\ge 0$:
\begin{align*}
\sup_{z\in (0,S(t))}|\partial_z m(t,z)|&\le M_1(t),\quad \sup_{z\in (0,S(t))} (-m(t,z)\partial_{z}^2 m(t,z))\le M_2(t).
\end{align*}
\item There exists a map $h:\R_{\ge 0}\times \R_{\ge 0}\to \R$ such that for fixed $t\ge 0$, $m(t,\cdot)$ is \emph{induced} by $h(t,\cdot)$, i.e.,
\begin{align*}
m(t,z)\partial_{z}^2 h(t,z)&=1\quad\text{for all }z\in (0,S(t)).
\end{align*}
The map $h$ is assumed to admit the following bound:

For fixed $t\ge 0$ and all $z\in [0,S(t))$, one has
\begin{align*}
|h(t,z)|&\le H_1(t)z^{\alpha_1(t)}+H_2(t)z^{\alpha_2(t)+1},
\end{align*}
where $H_1(t),H_2(t)\ge 0$, $\alpha_1(t)\in \left(\frac{d}{d+2},1\right)$ and $\alpha_2(t)\in [0,1]$ all depend continuously on $t$.\\

$h$ is locally Lipschitz w.r.t. $t$ in the following sense:

For all $T>0$, all $0\le t\le t'\le T$ and all $z\in [0,S(t)]$:
\begin{align*}
|h(t',z)-h(t,z)|\le \big(L_1(T)z^{\beta_1(T)}+L_2(T)z^{\beta_2(T)+1}\big)(t'-t),
\end{align*}
where $L_1(T),L_2(T)\ge 0$, $\beta_1(T)\in \left(\frac{d}{d+2},1\right)$ and $\beta_2(T)\in [0,1]$ all depend continuously on $T$.
\end{enumerate}
\end{definition}

The conditions (M1) and (M3) are concerned with the behaviour of $m$ for fixed time variable $t$: whereas (M1) is needed to be able to define a metric $\W_{m(t,\cdot)}$ on $\X(t)$, assumption (M3)---which is a time-dependent version of the so-called \emph{Lipschitz-semiconcavity (LSC)} condition from \cite{lisini2012} (see also \cite{zinsl2016})---yields $\lambda$-convexity along geodesics of the viscous regularization of the potential energy $u\mapsto \int_\Rd u\eta\dd x$ for some $\lambda\le 0$. This will be a cornerstone in our derivation of a time-discrete approximate weak formulation satisfied by the discrete solution $u_\tau$. In contrast, (M2) and (M4) are concerned with the behaviour of $m$ w.r.t. time $t$. There exists, at each $t\ge 0$ separately, one distinguished $0$-geodesically convex functional $\calH_t$,
\begin{align*}
\calH_t(u):=\int_\Rd h(t,u(x))\dd x,
\end{align*}
the so-called \emph{heat entropy} \cite{zm2014} since it induces the heat flow as its $\W_{m(t,\cdot)}$-gradient flow \cite{dns2009,lisini2010,lisini2012,zm2014}. Assumption (M4) is used to control the behaviour of $\calH_t$ with respect to time $t$, guaranteeing sufficient spatial regularity of the (time-)discrete solution $u_\tau$.\\

Our energy functionals are assumed to satisfy one of the following two conditions.\\

\begin{definition}[Admissible energy functionals]~
\begin{enumerate}[({E}1)]
\item \emph{Second order equations \eqref{eq:pde2o}:} Let $\ent$ be of the form \eqref{eq:ent2o},
where $f\in C^2(\R_{\ge 0})$ with $f(0)=0$, $f'(0)=0$ and $f''(z)\in [\gamma_0,\gamma_1]$ for $0<\gamma_0\le \gamma_1<\infty$ and all $z\ge 0$, and $\phi\in C^2(\Rd)$ is bounded from below and grows at most quadratically, with $\Delta \phi\in L^\infty(\Rd)$. 
\item \emph{Fourth order equations \eqref{eq:pde4o}:} Let $\ent$ be of the form \eqref{eq:ent4o}, 
with $\phi$ as in (E1) and $f\in C^2(\Rd\times \R_{\ge 0})$, $f(0,0)=0$, $\nabla_{(p,z)}f(0,0)=0$ and $\gamma_0\Id_{d+1}\le \nabla^2_{(p,z)} f(p,z)\le \gamma_1\Id_{d+1}$ for all $p\in\Rd$ and all $z\ge 0$, with $0<\gamma_0\le \gamma_1<\infty$.
\end{enumerate}
\end{definition}

We obtain the following on the convergence behaviour of $u_\tau$ as $\tau\searrow 0$.

\begin{thm}[Existence: second order case (E1)]\label{thm:exist2o}
Assume that the mobility satisfies \textup{(M1)--(M4)} and the energy functional $\ent$ is of the form \textup{(E1)} and let an initial datum $u_0\in \X\cap L^2(\Rd)$ with $u_0(x)\in [0,S(0)]$ for almost every $x\in\Rd$ be given. Then, for each $\tau>0$, the map $u_\tau$ obtained via \eqref{eq:mms} is well-defined. Furthermore, for each vanishing sequence $\tau_k\searrow 0$ $(k\to\infty)$, there exists a (nonrelabelled) subsequence and a limit map $u:\R_{\ge 0}\to\X$ such that the following is true for each fixed $T>0$:
\begin{enumerate}[(a)]
\item $u\in L^\infty([0,T];L^2(\Rd))\cap L^2([0,T];H^1(\Rd))\cap C^{1/2}([0,T];\W_2)$, and for almost every $t\ge 0$, one has $u(t,x)\le S(t)$ for almost all $x\in\Rd$;
\item $u_{\tau_k}(t)\to u(t)$ weakly$\ast$ in $\prb$, at each fixed $t>0$;
\item $u_{\tau_k}$ converges to $u$ strongly in $L^2([0,T];L^2(\Omega))$ for each bounded domain $\Omega\subset\Rd$ and weakly in \break $L^2([0,T];H^1(\Rd))$;
\item $u$ is a solution to \eqref{eq:pde2o} in the sense of distributions and $u(0,\cdot)=u_0$ a.e. on $\Rd$.
\end{enumerate}
\end{thm}

A similar statement also holds in the fourth order case.

\begin{thm}[Existence: fourth order case (E2)]\label{thm:exist4o}
Assume that the mobility satisfies \textup{(M1)--(M4)} and the energy functional $\ent$ is of the form \textup{(E2)} and let an initial datum $u_0\in \X\cap H^1(\Rd)$ with $u_0(x)\in [0,S(0)]$ for almost every $x\in\Rd$ be given. Then, for each $\tau>0$, the map $u_\tau$ obtained via \eqref{eq:mms} is well-defined. Furthermore, for each vanishing sequence $\tau_k\searrow 0$ $(k\to\infty)$, there exists a (nonrelabelled) subsequence and a limit map $u:\R_{\ge 0}\to\X$ such that the following is true for each fixed $T>0$:
\begin{enumerate}[(a)]
\item $u\in L^\infty([0,T];H^1(\Rd))\cap L^2([0,T];H^2(\Rd))\cap C^{1/2}([0,T];\W_2)$, and for almost every $t\ge 0$, one has $u(t,x)\le S(t)$ for almost all $x\in\Rd$;
\item $u_{\tau_k}(t)\to u(t)$ weakly$\ast$ in $\prb$, at each fixed $t>0$;
\item $u_{\tau_k}$ converges to $u$ strongly in $L^2([0,T];H^1(\Omega))$ for each bounded domain $\Omega\subset\Rd$ and weakly in $L^2([0,T];H^2(\Rd))$;
\item $u$ is a solution to \eqref{eq:pde4o} in the sense of distributions and $u(0,\cdot)=u_0$ a.e. on $\Rd$.
\end{enumerate}
\end{thm}

Note that in both cases, we do neither obtain uniqueness of solutions nor a monotonicity property of the energy functional in the limit $\tau\searrow 0$, since our notion of solution is very weak and the minimization problem in \eqref{eq:mms} lacks convexity. Notice furthermore that our results in principle also hold on a bounded and convex spatial domain. We refer to Section \ref{sec:gen} for a sketch of possible extensions in that direction.

\subsection{Strategy of proof and relation to the literature}
Our strategy of proof can be summarized as follows. Given the conditions (E1)/(E2) and (M1)--(M4), the map $u_\tau$ obtained via the modified minimizing movement scheme \eqref{eq:mms} is well-defined. It obeys a certain regularity property which is crucial for the passage to the limit as $\tau\searrow 0$ afterwards, since, due to the nonlinearity of the problem, weak convergence alone is not sufficient: there, we apply the so-called \emph{flow interchange} technique introduced in \cite{matthes2009}. Similarly as in \cite{lisini2012,zm2014}, the necessary auxiliary flow is given as the $\W_{m(n\tau,\cdot)}$-gradient flow of the \emph{(time-dependent) heat entropy} functional 
\begin{align*}
\calH_{n\tau}(u)=\int_\Rd h(n\tau,u)\dd x,
\end{align*}
where $h$ is the function from (M4). This auxiliary flow is---by construction---the heat flow. Nonautonomous evolution equations of gradient flow type have already been studied in \cite{ferreira2015,plazotta2016} from the opposite point of view: there, time-dependent energy functionals on time-independent metric spaces were considered and a different modification of the minimizing movement scheme was investigated. In \cite{plazotta2016}, where general non-convex problems have been studied, a certain Lipschitz condition for the free energy functional comes into play which resembles our additional conditions on $h$ from (M4). In \cite{ferreira2015,plazotta2016}, re-proving the classical properties of the minimizing movement scheme is more involved than for our scheme from Definition \ref{def:minmov}---in contrast, deriving higher regularity estimates brings additional difficulties in our case, since now the heat entropy is a time-dependent functional. Non-autonomous equations of Wasserstein gradient flow form with linear mobility have also been investigated in \cite{petrelli2004} using a time-averaged form of the classical minimizing movement scheme. Notice that our equations do not in general allow for a time-dependent scaling which transforms the problem into autonomous form (for studies in this direction, see e.g. \cite{bdik2007}).

Nonautonomous linear and semilinear equations equations have, in contrast, been investigated more exhaustively using semigroup theory (see e.g. \cite{aulbach1996} and references therein). More recently, properties of time-dependent Riemannian manifolds have been studied \cite{sturm2016}, and corresponding evolution problems have been investigated \cite{lengeler2013}.

The partial Riemannian structure on the space $\X(t)$ induced by the distance $\W_{m(t,\cdot)}$ was introduced in \cite{dns2009,lisini2010} and later generalized to the vector-valued framework in \cite{zm2014}; the structure of geodesics on that space was first investigated in \cite{carrillo2010}. In \cite{lisini2012}, the formal gradient flow with respect to that distance was used to prove existence of solutions for a certain class of Cahn-Hilliard type fourth-order equations. There, also mobility functions which are not Lipschitz continuous w.r.t. $z$ could be used using an approximation technique. In this work (see also \cite{zinsl2016}), we will employ this strategy to extend our results from Section \ref{subsec:main} to wider classes of mobility functions as those satisfying (M1)--(M4), see Section \ref{sec:gen}.

\subsection{Plan of the paper}
The paper is organized as follows. In Section \ref{sec:pre}, we summarize some preliminary facts on gradient flows with respect to the $L^2$-Wasserstein distance and its modifications. Afterwards, the variational scheme from Definition \ref{def:minmov} is studied in Section \ref{sec:minmov}. Section \ref{sec:conv} then is concerned with the derivation of the approximate discrete weak formulation and the passage to the continuous-time limit $\tau\searrow 0$. We sketch several possible generalizations in Section \ref{sec:gen}.


\section{Preliminaries}\label{sec:pre}
Derivatives with respect to state variables $z$ or $p$ will always be indicated as $\partial_z$ or $\nabla_p$, respectively. Spatial derivatives $\nabla,~\dv$ or $\Delta$ are denoted without further indication. We now briefly summarize some results from the theory of gradient flows in spaces of measures. Our ambient space $\X$ from \eqref{eq:space} is---by the slight abuse of notation where we identify probability densities $u$ with their corresponding probability measures $u\cdot\mathscr{L}^d$---a subspace of the space of probability measures $\prb$ on $\Rd$. A sequence $(\mu_k)_{k\in\N}$ in $\prb$ is said to converge \emph{weakly}$\ast$ to some $\mu\in\prb$ if for all continuous and bounded maps $g:\Rd\to\R$, one has
\begin{align*}
\lim_{k\to\infty}\int_\Rd g\dd\mu_k&=\int_\Rd g\dd\mu.
\end{align*}
The space $\prb$ can be endowed with the so-called $L^2$-Wasserstein (pseudo-)distance: For each $\mu_0,\mu_1\in\prb$, 
\begin{align*}
\W_2(\mu_0,\mu_1)=\inf\left\{\int_{\Rd\times\Rd}|x-y|\dd\gamma(x,y):~\gamma\in\Gamma(\mu_0,\mu_1)\right\}^{1/2},
\end{align*}
where $\Gamma(\mu_0,\mu_1)$ denotes the set of transport plans from $\mu_0$ to $\mu_1$ (for more details on optimal transport, see \cite{villani2003}). A dynamical characterization of $\W_2$---which will be made use of here---was found by Benamou and Brenier \cite{brenier2000}:
\begin{align*}
&\W_{2}(u_0,u_1)=\inf\left\{\int_0^1\int_\Rd \frac{|w_s|^2}{u_s}\dd x \dd s:~(u_s,w_s)_{s\in[0,1]}\in \mathscr{C},~u_s|_{s=0}=u_0,~u_s|_{s=1}=u_1\right\}^{1/2},
\end{align*}
where $\mathscr{C}$ is a suitable subclass of solutions to the \emph{continuity equation} $\partial_s u_s=-\dv w_s$ in the sense of distributions on $[0,1]\times\Rd$. This characterization has been the starting point for the definition of the distances $\W_{m(t,\cdot)}$ in \eqref{eq:metr} for \emph{nonlinear} mobility functions \cite{dns2009,lisini2010}.\\

In order to derive the necessary \emph{a priori} estimates on the discrete solution $u_\tau$ obtained via our variational scheme from Definition \ref{def:minmov}, we use the so-called \emph{flow interchange} lemma from \cite{matthes2009}. In advance of its precise statement, we introduce the following notion of gradient flow (compare \cite{savare2008,daneri2008}).
\begin{definition}[$\kappa$-flows]
Let $\auxil:\prb\to\R\cup\{+\infty\}$ be a proper and lower semicontinuous functional on the (pseudo-)metric space $(\prb,\metr)$ and let $\kappa\in\R$. A continuous semigroup $\flow^{\auxil}$ on $(\prb,\metr)$ satisfying
the \emph{evolution variational estimate}
 \begin{align*}
   \frac{1}{2}\frac{\dd^+}{\dd s}\metr^2(\flow_{s}^{\auxil}(w),\tilde w)+\frac{\kappa}{2}\metr^2(\flow_{s}^{\auxil}(w),\tilde w)+\auxil(\flow_{s}^{\auxil}(w))
   &\le\auxil(\tilde w)
\end{align*}
for arbitrary $w,\tilde w\in\operatorname{Dom}(\mathcal{A})$ and for all $s\ge 0$, as well as the monotonicity condition
\begin{align*}
\auxil(\flow_t^\auxil(w))\le \auxil(\flow_s^\auxil(w))\quad\forall 0\le s\le t
\end{align*}
for all $w\in\prb$,
is called \emph{$\kappa$-flow} or \emph{gradient flow} of $\auxil$.
\end{definition}
Notice that \cite{daneri2008,liero2013} if $\auxil$ induces a $\kappa$-flow on $(\prb,\metr)$ for some $\kappa\in\R$, then $\auxil$ is $\kappa$-convex along geodesics in $(\prb,\metr)$ in the sense of McCann \cite{mccann1997}.
\begin{thm}[Flow interchange lemma {\cite[Thm. 3.2]{matthes2009}}]\label{thm:flowinterchange}
Let $\mathcal{B}$ be a proper and lower semicontinuous functional on $(\prb,\metr)$ and assume that there exists a $\lambda$-flow $\flow^\mathcal{B}$ for some $\lambda\in\R$. Let furthermore $\mathcal{A}$ be another proper, lower semicontinuous functional on $(\prb,\metr)$ such that $\operatorname{Dom}(\mathcal{A})\subset \operatorname{Dom}(\mathcal{B})$. Assume that, if $\tau>0$ and $\tilde w\in \prb$ are such that the functional $u\mapsto \frac1{2\tau}\metr^2(u,\tilde w)+\auxil(u)$ is proper, it also possesses a minimizer $w$ on $\prb$.
Then, the following holds:
\begin{align*}
\mathcal{B}(w)+\tau \mathrm{D}^\mathcal{B}\mathcal{A}(w)+\frac{\lambda}{2}\metr^2(w,\tilde w)&\le \mathcal{B}(\tilde w).
\end{align*}
There, $\mathrm{D}^\mathcal{B}\mathcal{A}(w)$ denotes the \emph{dissipation} of the functional $\mathcal{A}$ along the $\lambda$-flow $\flow^{\mathcal{B}}$ of the functional $\mathcal{B}$, i.e.
\begin{align*}
\mathrm{D}^\mathcal{B}\mathcal{A}(w):=\limsup_{h\searrow 0}\frac{\mathcal{A}(w)-\mathcal{A}(\flow_{h}^{\mathcal{B}}(w)    )   }{h}.
\end{align*}
\end{thm}
With the uniform estimates on $u_\tau$ derived with the flow interchange lemma, one seeks to pass to the limit $\tau\searrow 0$ w.r.t. a suitably strong notion of convergence. To this end, the following extension of the classical Aubin-Lions compactness lemma to the metric setting is useful.
\begin{thm}[Extension of the Aubin-Lions lemma {\cite[Thm. 2]{rossi2003}}]\label{thm:ex_aub}
Let $\ban$ be a Banach space and let \break $\mathcal{A}:\,\ban\to[0,\infty]$ be lower semicontinuous and have relatively compact sublevels in $\ban$. Let furthermore \break $\W:\,\ban\times \ban\to[0,\infty]$ be lower semicontinuous and such that $\W(u,\tilde u)=0$ for $u,\tilde u\in \dom(\mathcal{A})$ implies $u=\tilde u$.

If for a sequence $(U_k)_{k\in\N}$ of measurable functions $U_k:\,(0,T)\to \ban$, one has
\begin{align}
\sup_{k\in\N}\int_0^T\mathcal{A}(U_k(t))\dd t&<\infty\qquad\qquad\qquad\text{and}\label{eq:hypo1}\\
\lim_{h\searrow 0}\sup_{k\in\N}\int_0^{T-h}\W(U_k(t+h),U_k(t))\dd t&=0,\label{eq:hypo2}
\end{align}
then there exists a subsequence that converges in measure w.r.t. $t\in(0,T)$ to a limit~$U:\,(0,T)\to \ban$.
\end{thm}
We conclude this preliminary section with some elementary properties of the mobility functions $m$ and their induced distances.
\begin{prop}[Properties of $m$ and $\W_{m(t,\cdot)}$]\label{prop:propm}
The following statements hold:
\begin{enumerate}[(a)]
\item For each fixed $t\ge 0$, the distance functional $\W_{m(t,\cdot)}$ is in both arguments lower semicontinuous with respect to weak$\ast$-convergence in $\prb$.
\item For all $0\le t_0\le t_1$, one has
\begin{align*}
\X(t_0)\subset \X(t_1)\subset \X.
\end{align*}
\item For all $T>0$, all $t\in [0,T]$ and all $z\in [0,S(t))$, one has, for $C(T):=\max\limits_{t\in[0,T]}\lim\limits_{z'\searrow 0}\partial_z m(t,z')<\infty$:
\begin{align*}
\partial_z m(t,z)&\le C(T),\\
m(t,z)&\le C(T)z=:\bar m(T,z).
\end{align*}
\item At each $T>0$, $\W_{\bar m(T,\cdot)}$ is a scalar multiple of the $L^2$-Wasserstein distance on the space of probability measures and the following estimate holds for $t\in[0,T]$ and $u,\tilde u\in\X(t)$:
\begin{align*}
\W_2(u,\tilde u)\le \sqrt{C(T)}\W_{\bar m(T,\cdot)}(u,\tilde u)\le \sqrt{C(T)}\W_{m(t,\cdot)}(u,\tilde u).
\end{align*}
\end{enumerate}
\end{prop}
\begin{proof}
The lower semicontinuity of $\W_{m(t,\cdot)}$ is classical, see \cite{dns2009,lisini2010}. Obviously, (b) holds as $S$ is nondecreasing (M1). The estimates in (c) are a straightforward consequence of (M1)--(M3), since $\partial_z m(t,\cdot)$ is nonincreasing. Finally, (d) can easily be derived using (c) and the characterization of the $L^2$-Wasserstein distance by the Benamou-Brenier formula.
\end{proof}


\section{The variational scheme}\label{sec:minmov}
In this section, we prove the well-posedness of the scheme from Definition \ref{def:minmov}, i.e., the sequence of successive minimizers $(u_\tau^n)_{n\in\N}$ in \eqref{eq:mms} exists, and investigate certain regularity properties. First, some elementary properties of the free energy $\ent$ are summarized.
\begin{lemma}[Properties of $\ent$, case (E1)]\label{eq:propent1}
The following statements hold:
\begin{enumerate}[(a)]
\item For all $u\in\X\cap L^2(\Rd)$, one has
\begin{align*}
-\infty<\inf_{x\in\Rd}\phi\le\frac{\gamma_0}{2}\|u\|_{L^2}^2+\inf_{x\in\Rd}\phi&\le \ent(u)\le \frac{\gamma_1}{2}\|u\|_{L^2}^2+C_\phi(\mom{u}+1)<\infty,
\end{align*}
where $C_\phi>0$ is such that $\phi(x)\le C_\phi(|x|^2+1)$ for all $x\in\Rd$.
\item Let a sequence $(u_k)_{k\in\N}$ in $\X\cap L^2(\Rd)$ be given and assume that $u_k$ converges weakly$\ast$ in $\prb$ and weakly in the space $L^2(\Rd)$ to some $u\in\X\cap L^2(\Rd)$. Then, 
\begin{align*}
\ent(u)&\le \liminf_{k\to\infty}\ent(u_k)<\infty.
\end{align*} 
\end{enumerate}
\end{lemma}
\begin{proof}
By Taylor's theorem, we have $f(z)=f(z)+zf'(0)+\frac12 f''(\tilde z)z^2$ for all $z$ and some $\tilde z\in (0,z)$. Hence, using the bounds on $f''$ from (E1) and the estimate $-\infty<\inf\phi\le \phi(x)\le C_\phi(|x|^2+1)<\infty$, the claim in (a) follows immediately. Notice in particular that $f$ is a nonnegative and convex function. Therefore, the map $u\mapsto \int_\Rd f(u(x))\dd x$ is weakly lower semicontinuous in $L^2(\Rd)$. Together with the---as $\phi$ is bounded from below---obvious weak$\ast$-lower semicontinuity of $u\mapsto \int_\Rd \phi(x)u(x)\dd x$ \cite[Lemma 5.1.7]{savare2008}, part (b) follows.
\end{proof}
Similar arguments show an analogous result in the case of gradient-dependent energy of the form in (E2):
\begin{lemma}[Properties of $\ent$, case (E2)]\label{eq:propent2}
The following statements hold:
\begin{enumerate}[(a)]
\item For all $u\in\X\cap H^1(\Rd)$, one has
\begin{align*}
-\infty<\inf_{x\in\Rd}\phi\le\frac{\gamma_0}{2}\|u\|_{H^1}^2+\inf_{x\in\Rd}\phi&\le \ent(u)\le \frac{\gamma_1}{2}\|u\|_{H^1}^2+C_\phi(\mom{u}+1)<\infty,
\end{align*}
where $C_\phi>0$ is such that $\phi(x)\le C_\phi(|x|^2+1)$ for all $x\in\Rd$.
\item Let a sequence $(u_k)_{k\in\N}$ in $\X\cap H^1(\Rd)$ with uniformly bounded second moments $\mom{u_k}$ be given and assume that $u_k$ converges weakly$\ast$ in $\prb$ and weakly in the space $H^1(\Rd)$ to some $u\in\X\cap H^1(\Rd)$. Then, 
\begin{align*}
\ent(u)&\le \liminf_{k\to\infty}\ent(u_k)<\infty.
\end{align*} 
\end{enumerate}
\end{lemma}
\subsection{Well-posedness and classical properties}
This paragraph is concerned with the well- \break posedness of the scheme \eqref{eq:mms} and properties resembling the classical estimates known for the autonomous case. 
\begin{prop}[Well-posedness and classical estimates, case (E1)]\label{prop:mmsclass1}
Assume that $u_0\in \X(0)\cap L^2(\Rd)$ and let $\bar\tau>0$. Then, for all $\tau\in (0,\bar\tau]$ and all $n\in\N$, the map
\begin{align}
\ent_n(u):=\frac1{2\tau}\W_{m(n\tau,\cdot)}(u_\tau^{n-1},u)^2+\ent(u)\label{eq:yosida}
\end{align}
possesses a minimizer $u_\tau^n$ on $\X$ belonging to $\X(n\tau)\cap L^2(\Rd)$. Furthermore, the following estimates hold:
\begin{align}
\ent(u_\tau^n)&\le \ent(u_\tau^{n-1})\le \ent(u_0)<\infty,\quad\text{for all }n\in\N,\label{eq:enest}\\
\sum_{n=1}^\infty \W_{m(n\tau,\cdot)}(u_\tau^n,u_\tau^{n-1})^2&\le 2\tau (\ent(u_0)-\inf_{u\in\X}\ent(u))<\infty,\label{eq:distest}\\
\W_2(u_\tau(t),u_\tau(s))&\le \sqrt{2C(T+\bar\tau)(\ent(u_0)-\inf_{u\in\X}\ent(u))\max(\tau,|t-s|)},\label{eq:holdest}\\&\quad\text{for all }s,t\in [0,T]\text{ and all }T>0.\nonumber
\end{align}
\end{prop}
\begin{proof}
We proceed by induction on $n\in\N$ and consider a minimizing sequence $(u_k)_{k\in\N}$ for $\ent_n$ which consequently is a sequence in $\X(n\tau)\cap L^2(\Rd)$. Using the bounds from Proposition \ref{prop:propm}(b), one has for a suitable constant $C>0$ that $\|u_k\|_{L^2}\le C$ and $\mom{u_k}\le C$ for all $k\in\N$. Subsequently, the Banach-Alaoglu and Prokhorov theorems yield the existence of a (non-relabelled) subsequence and a limit $u\in \X(n\tau)\cap L^2(\Rd)$ such that $u_k$ converges to $u$ both weakly$\ast$ in $\prb$ and weakly in $L^2(\Rd)$. Since $\ent$ and $\W_{m(n\tau,\cdot)}(u_\tau^{k-1},\cdot)$ are lower semicontinuous w.r.t. these convergences, the limit $u$ indeed is a minimizer of $\ent_n$.

The minimizing property of $u_\tau^n$ immediately yields \eqref{eq:enest}, and after summing up, also \eqref{eq:distest}. For the H\"older type estimate \eqref{eq:holdest}, we fix $T>0$ and $s,t\in [0,T]$ and find, using the triangle inequality and Proposition \ref{prop:propm}(d):
\begin{align*}
\W_2(u_\tau(t),u_\tau(s))&\le \sum_{n=M}^N \W_2(u_\tau^n,u_\tau^{n-1})\le \sqrt{C(T+\bar\tau)}\sum_{n=M}^N \W_{m(n\tau,\cdot)}(u_\tau^n,u_\tau^{n-1})
\end{align*}
for some $M,N\in\N$ with $M\le N$. We proceed using H\"older's inequality and \eqref{eq:distest} to obtain the desired estimate:
\begin{align*}
\sqrt{C(T+\bar\tau)}\sum_{n=M}^N \W_{m(n\tau,\cdot)}(u_\tau^n,u_\tau^{n-1})&\le \left(C(T+\bar\tau)(N-M+1)\sum_{n=M}^N\W_{m(n\tau,\cdot)}(u_\tau^n,u_\tau^{n-1})^2\right)^{1/2}\\
&\le \sqrt{2C(T+\bar\tau)(\ent(u_0)-\inf_{u\in\X}\ent(u))\max(\tau,|t-s|)}.\qedhere
\end{align*}
\end{proof}
By an easy adaptation of the proof, one obtains
\begin{prop}[Well-posedness and classical estimates, case (E2)]\label{prop:mmsclass2}
Assume that $u_0\in \X(0)\cap H^1(\Rd)$ and let $\bar\tau>0$. Then, for all $\tau\in (0,\bar\tau]$ and all $n\in\N$, the map from \eqref{eq:yosida}
possesses a minimizer $u_\tau^n$ on $\X$ belonging to $\X(n\tau)\cap H^1(\Rd)$, and the estimates \eqref{eq:enest}--\eqref{eq:holdest} hold.
\end{prop}

\subsection{Additional regularity}
The minimizers of $\ent_n$ enjoy a certain regularity property which is crucial for the passage to the continuous-time limit afterwards. To this end, we introduce for $t\ge 0$ the \emph{time-dependent heat entropy}
\begin{align}\label{eq:heatent}
\calH_{(\cdot)}:\R_{\ge 0}\times\X\to\R_\infty,\quad \calH_t(u):=\begin{cases}\int_\Rd h(t,u(t,x))\dd x&\text{if }u\in \X(t),\\ +\infty&\text{otherwise.}\end{cases}
\end{align}
We first prove some elementary properties of $\calH$:
\begin{lemma}[Time-dependent heat entropy]\label{lem:heatent}
The following statements hold :
\begin{enumerate}[(a)]
\item There exists a continuous function $C:\R_{\ge 0}\to\R_{\ge 0}$ such that for each fixed $t\ge 0$ and all $u\in\X(t)\cap L^2(\Rd)$:
\begin{align}\label{eq:carl}
|\calH_t(u)|&\le C(t)(\|u\|_{L^2}^2+\mom{u}+1)<\infty.
\end{align}
\item For each fixed $t\ge 0$, $\calH_t$ induces a $0$-flow $\flow^{\calH_{t}}$ on $\X(t)$ coinciding with the heat flow, viz.
\begin{align*}
\partial_s\flow^{\calH_{t}}_s(u)&=\Delta \flow^{\calH_{t}}_s(u),\quad \flow^{\calH_{t}}_0(u)=u,\quad\text{for }u\in\X(t).
\end{align*}
\end{enumerate}
\end{lemma}
\begin{proof}
Thanks to the bounds on $h$ from (M4), one has for all $t\ge 0$ and all $u\in\X(t)$ that
\begin{align*}
\calH_t(u)&\le H_1(t)\int_\Rd u^{\alpha_1(t)}\dd x+H_2(t)\int_\Rd u^{\alpha_2(t)+1}\dd x.
\end{align*}
The first integral can be estimated using H\"older's inequality (recall that $\alpha_1(t)<1$):
\begin{align*}
\int_\Rd u^{\alpha_1(t)}\dd x&=\int_\Rd (u(x)(|x|^2+1))^{\alpha_1(t)}(|x|^2+1)^{-\alpha_1(t)}\dd x\\&\le (\mom{u}+1)^{\alpha_1(t)}\left(\int_\Rd (|x|^2+1)^{-\frac{\alpha_1(t)}{1-\alpha_1(t)}}\dd x\right)^{1-\alpha_1(t)},
\end{align*}
and $\int_\Rd (|x|^2+1)^{-\frac{\alpha_1(t)}{1-\alpha_1(t)}}\dd x$ is finite since $\frac{-2\alpha_1(t)}{1-\alpha_1(t)}+d-1<-1$ thanks to $\alpha_1(t)>\frac{d}{d+2}$. The second integral above can be controlled with the $L^p$ interpolation inequality as follows:
\begin{align*}
\int_\Rd u^{\alpha_2(t)+1}\dd x\le \|u\|_{L^1}^{1-\theta(t)}\|u\|_{L^2}^{\theta(t)}=\|u\|_{L^2}^{\theta(t)},
\end{align*}
for a suitable $\theta(t)\in(0,1)$ depending continuously on $t$. All in all, applying Young's inequality provides \eqref{eq:carl}. 

Since $m(t,\cdot)$ is induced by $h(t,\cdot)$ and since $\calH_t$ is finite on $\X(t)$, a $0$-flow on $\X(t)$ is induced by $\calH_t$, for each fixed $t\ge 0$ (see \cite{lisini2012,zm2014}). By construction, $\flow^{\calH_t}$ is the heat flow. 
\end{proof}
\begin{prop}[Additional regularity, case (E1)]\label{prop:addreg1}
Assume that $u_0\in \X(0)\cap L^2(\Rd)$ and let $\bar\tau>0$. Then, there exists $C>0$ such that for all $\tau\in (0,\bar\tau]$ and all $n\in\N$, one has
\begin{align}\label{eq:addreg1}
\tau \|\nabla u_\tau^n\|_{L^2}^2\le C(\calH_{n\tau}(u_\tau^{n-1})-\calH_{n\tau}(u_\tau^{n})+\tau).
\end{align}
\end{prop}
\begin{proof}
Our proof is a (by now almost classical) application of the \emph{flow interchange} technique from \cite{matthes2009}. The suitable auxiliary flow is the aforementioned $0$-flow $\flow^{\calH_{n\tau}}$ of the functional $\calH_t$ at $t=n\tau$, for each fixed $n\in\N$. We calculate the dissipation of $\ent$ along the semigroup $(\flow^{\calH_{n\tau}}_s(u_\tau^n))_{s\ge 0}$ and write $u_s:=\flow^{\calH_{n\tau}}_s(u_\tau^n)$ for brevity:
\begin{align*}
-\frac{\dd}{\dd s}\ent(u_s)=-\int_\Rd(f'(u_s)+\phi)\Delta u_s\dd x.
\end{align*}
Integrating by parts and using (E1), we obtain
\begin{align*}
-\frac{\dd}{\dd s}\ent(u_s)&=\int_\Rd(f''(u_s)|\nabla u_s|^2-u_s\Delta\phi )\dd x\ge \gamma_0\int_\Rd |\nabla u_s|^2\dd x-\|\Delta\phi\|_{L^\infty}.
\end{align*}
As $s\searrow 0$, one has
\begin{align*}
\dff^{\calH_{n\tau}}\ent(u_\tau^n)&\ge\liminf_{s\searrow 0}\left(-\frac{\dd}{\dd s}\ent(u_s)\right)\ge \gamma_0\int_\Rd |\nabla u_\tau^n|^2\dd x-\|\Delta\phi\|_{L^\infty}.
\end{align*}
We apply the flow interchange lemma (Theorem \ref{thm:flowinterchange}) and rearrange to obtain the desired result:
\begin{align*}
\calH_{n\tau}(u_\tau^n)+\tau\left[\gamma_0\int_\Rd |\nabla u_s|^2\dd x-\|\Delta\phi\|_{L^\infty}\right]&\le \calH_{n\tau}(u_\tau^{n-1}). \qedhere
\end{align*}
\end{proof}
A similar result is also true for gradient-dependent energy:
\begin{prop}[Additional regularity, case (E2)]\label{prop:addreg2}
Assume that $u_0\in \X(0)\cap H^1(\Rd)$ and let $\bar\tau>0$. Then, there exists $C>0$ such that for all $\tau\in (0,\bar\tau]$ and all $n\in\N$, one has
\begin{align}\label{eq:addreg2}
\tau \|\nabla^2 u_\tau^n\|_{L^2}^2\le C(\calH_{n\tau}(u_\tau^{n-1})-\calH_{n\tau}(u_\tau^{n})+\tau).
\end{align}
\end{prop}
\begin{proof}
As in the proof of Proposition \ref{prop:addreg1}, we derive the dissipation of $\ent$ along $(\flow^{\calH_{n\tau}}_s(u_\tau^n))_{s\ge 0}$, integrate by parts and use (E2):
\begin{align*}
-\frac{\dd}{\dd s}\ent(u_s)&=-\int_\Rd\left[-\dv\nabla _pf(\nabla u_s,u_s)+\partial_z f(\nabla u_s,u_s)+\phi\right]\Delta u_s\dd x\\
&=-\int_\Rd\left[\sum_{i=1}^d(-\dv\nabla_pf(\nabla u_s,u_s)+\partial_z f(\nabla u_s,u_s))\partial_{x_i}\partial_{x_i}u_s+u_s\Delta\phi\right]\dd x\\
&=-\int_\Rd\left[\sum_{i=1}^d\nabla_p f(\nabla u_s,u_s)\cdot\partial_{x_i}\partial_{x_i}\nabla u_s+\partial_z f(\nabla u_s,u_s)\partial_{x_i}\partial_{x_i}u_s+u_s\Delta\phi\right]\dd x\\
&=\int_\Rd \left[\sum_{i=1}^d\partial_{x_i}\nabla_pf(\nabla u_s,u_s)\cdot\partial_{x_i}\nabla u_s+\partial_{x_i}\partial_zf(\nabla u_s,u_s)\partial_{x_i}u_s-u_s\Delta\phi\right]\dd x\\
&=\int_\Rd \left[\sum_{i=1}^d\begin{pmatrix}\partial_{x_i}\nabla u_s \\ \partial_{x_i}u_s\end{pmatrix}^\tT \nabla_{(p,z)}^2f(\nabla u_s,u_s)\begin{pmatrix}\partial_{x_i}\nabla u_s \\ \partial_{x_i}u_s\end{pmatrix}-u_s\Delta\phi\right]\dd x\\
&\ge \int_\Rd\left[\gamma_0\sum_{i=1}^d((\partial_{x_i}u_s)^2+|\partial_{x_i}\nabla u_s|^2)-u_s\Delta\phi\right]\dd x\\
&\ge \gamma_0 \|\nabla^2 u_s\|_{L^2}^2-\|\Delta\phi\|_{L^\infty}.
\end{align*}
From there, the proof goes along the same lines as the proof of Proposition \ref{prop:addreg1}.
\end{proof}
We can now summarize the relevant \emph{a priori} estimates on the discrete solution $u_\tau$.
\begin{lemma}[\emph{A priori} estimates, case (E1)]\label{lem:apri1}
Let $u_0\in \X(0)\cap L^2(\Rd)$, $\bar\tau>0$ and $T>0$ be given. Then, there exists a constant $C>0$ such that for all $\tau\in(0,\bar\tau]$:
\begin{enumerate}[(a)]
\item $\|u_\tau\|_{L^\infty([0,T];L^2)}\le C,$
\item $\sup\limits_{t\in[0,T]}\mom{u_\tau(t,\cdot)}\le C,$
\item $\|u_\tau\|_{L^2([0,T];H^1)}\le C.$
\end{enumerate}
\end{lemma}
\begin{proof}
The parts (a) and (b) are direct consequences of the estimates \eqref{eq:enest} and \eqref{eq:holdest}, respectively. For part (c), we use \eqref{eq:addreg1} to see for $N:=\lfloor\frac{T}{\tau}\rfloor+1$ that
\begin{align}\label{eq:apriori3_1}
\|\nabla u_\tau\|_{L^2([0,T];L^2)}^2&\le\tau\sum_{n=1}^N \|\nabla u_\tau^n\|_{L^2}^2\le C\sum_{n=1}^N (\calH_{n\tau}(u_\tau^{n-1})-\calH_{n\tau}(u_\tau^{n})+\tau).
\end{align}
Inserting a suitable term and employing the local Lipschitz condition for $h$ from (M4) yields for a fixed $n$:
\begin{align*}
&\calH_{n\tau}(u_\tau^{n-1})-\calH_{n\tau}(u_\tau^{n})=\calH_{n\tau}(u_\tau^{n-1})-\calH_{(n-1)\tau}(u_\tau^{n-1})+\calH_{(n-1)\tau}(u_\tau^{n-1})-\calH_{n\tau}(u_\tau^{n})\\
&\quad\le \tau\int_\Rd\left[L_1(T+\bar\tau)(u_\tau^{n-1})^{\beta_1(T+\bar\tau)}+L_2(T+\bar\tau)(u_\tau^{n-1})^{\beta_2(T+\bar\tau)+1}\right]\dd x+\calH_{(n-1)\tau}(u_\tau^{n-1})-\calH_{n\tau}(u_\tau^{n}).
\end{align*}
Thus, we get from \eqref{eq:apriori3_1} simplifying the telescopic sum that
\begin{align*}
\|\nabla u_\tau\|_{L^2([0,T];L^2)}^2&\le \tau\sum_{n=1}^N C\left[1+\int_\Rd\left[L_1(T+\bar\tau)(u_\tau^{n-1})^{\beta_1(T+\bar\tau)}+L_2(T+\bar\tau)(u_\tau^{n-1})^{\beta_2(T+\bar\tau)+1}\right]\dd x\right]\\
&\quad+C(\calH_0(u_0)-\calH_{N\tau}(u_\tau^N)).
\end{align*}
Proceeding as in the proof of \eqref{eq:carl} from Lemma \ref{lem:heatent} gives us for a constant $C'>0$ that
\begin{align*}
&\|\nabla u_\tau\|_{L^2([0,T];L^2)}^2\le C(T+\bar\tau)+C'\tau\sum_{n=0}^N \left[\mom{u_\tau^{n}}+\|u_\tau^{n}\|_{L^2}+1\right],
\end{align*}
which can be controlled by a finite constant independent of $\tau$ with the parts (a) and (b) of this proposition.
\end{proof}
Analogously, one shows:
\begin{lemma}[\emph{A priori} estimates, case (E2)]\label{lem:apri2}
Let $u_0\in \X(0)\cap H^1(\Rd)$, $\bar\tau>0$ and $T>0$ be given. Then, there exists a constant $C>0$ such that for all $\tau\in(0,\bar\tau]$:
\begin{enumerate}[(a)]
\item $\|u_\tau\|_{L^\infty([0,T];H^1)}\le C,$
\item $\sup\limits_{t\in[0,T]}\mom{u_\tau(t,\cdot)}\le C,$
\item $\|u_\tau\|_{L^2([0,T];H^2)}\le C.$
\end{enumerate}
\end{lemma}

\section{Weak formulation and passage to the continuous-time limit}\label{sec:conv}
First, we derive an approximate weak formulation of equation \eqref{eq:pde} corresponding to the Euler-Lagrange equation associated to the minimization problem in \eqref{eq:mms}. Afterwards, we show that $u_\tau$ converges in a suitable sense to a map $u$ which fulfills the (time-continuous) distributional formulation of equation \eqref{eq:pde}.
\subsection{Discrete weak formulation}
In this paragraph, we derive an approximate weak formulation for equation \eqref{eq:pde} in discrete time. Again, our method relies on the flow interchange lemma (Theorem \ref{thm:flowinterchange}). Here, the auxiliary flow is the $\lambda$-flow of the time-dependent \emph{regularized potential energy}
\begin{align*}
\aV_{(\cdot)}:\R_{\ge 0}\times\X\to\R_\infty,\quad \aV_t(u)=\sqrt{\tau}\calH_t(u)+\int_\Rd u\eta \dd x,
\end{align*}
for some test function $\eta\in C^\infty_c(\Rd)$ and some $\tau>0$:
\begin{lemma}[$t$-uniform $\lambda$-convexity of $\aV_t$]\label{lem:potcvx}
Let $T>0$, $\eta\in C^\infty_c(\Rd)$, $\bar\tau>0$ and $\tau\in(0,\bar\tau]$. There exists a constant $C>0$ such that for every $t\in [0,T]$ and with $\lambda:=-\frac{C}{\sqrt{\tau}}$, the functional $\aV_t$ is $\lambda$-convex and induces a $\lambda$-flow $\flow_s^{\aV_t}$ on $\X$. The flow $\flow_s^{\aV_t}$ solves the \emph{viscous continuity equation}, viz.
\begin{align*}
\partial_s \flow_s^{\aV_t}(u)=\sqrt{\tau}\Delta \flow_s^{\aV_t}(u)+\dv(m(t,\flow_s^{\aV_t}(u))\nabla\eta),\quad \flow_0^{\aV_t}(u)=u.
\end{align*}
\end{lemma}
\begin{proof}
According to the criterion from \cite{liero2013} (see also \cite{lisini2012,zm2014}), one has to show that for each fixed $t\in[0,T]$, all $z\in (0,S(t))$, $\zeta\in\R$ and all $w_1,w_2\in\R$ with $|w_1|\le \|\eta\|_{C^2}$, $|w_2|\le \|\eta\|_{C^2}$, the following holds for some $\lambda\in\R$:
\begin{align}\label{eq:ccrit}
-\frac{\sqrt{\tau}}{2}\partial_{z}^2m(t,z)\zeta^2+\frac12\partial_{z}^2m(t,z)m(t,z)w_1\zeta+\partial_{z}m(t,z)m(t,z)w_2&\ge \lambda m(t,z).
\end{align}
Using Young's inequality, we get
\begin{align*}
&-\frac{\sqrt{\tau}}{2}\partial_{z}^2m(t,z)\zeta^2+\frac12\partial_{z}^2m(t,z)m(t,z)w_1\zeta+\partial_{z}m(t,z)m(t,z)w_2\\
&\ge \left(\frac1{8\sqrt{\tau}}\partial_z^2 m(t,z)m(t,z)w_1^2+\partial_z m(t,z)w_2\right)m(t,z)\\
&\ge -m(t,z)\left(\frac1{8\sqrt{\tau}}M_2(t)w_1^2+M_1(t)w_2\right),
\end{align*}
with condition (M3) used in the last step. The bounds on $w_1$, $w_2$ and $\tau$ and the continuity of $M_1$ and $M_2$ yield
\begin{align*}
-m(t,z)\left(\frac1{8\sqrt{\tau}}M_2(t)w_1^2+M_1(t)w_2\right)\ge -\frac{m(t,z)}{\sqrt{\tau}}\left(\frac18 R^2 \max_{t\in[0,T]}M_2(t)+\sqrt{\bar\tau}R\max_{t\in[0,T]}M_1(t)\right),
\end{align*}
so \eqref{eq:ccrit} holds with $\lambda:=-\frac18 R^2 \max\limits_{t\in[0,T]}M_2(t)-\sqrt{\bar\tau}R\max\limits_{t\in[0,T]}M_1(t)$.
\end{proof}
\begin{lemma}[Discrete weak formulation, case (E1)]\label{lem:dweak1}
Let $u_0\in \X(0)\cap L^2(\Rd)$, $\bar\tau>0$ and let test functions $\eta\in C^\infty_c(\Rd)$ and $\psi\in C^\infty_c(\Rd)\cap C(\R_{\ge 0})$ be given. Then, there exists a constant $C>0$ such that for all $\tau\in (0,\bar\tau]$, one has
\begin{align}
\label{eq:dweak1}
\begin{split}
&\left|\int_0^\infty\int_\Rd \left[\frac{\psi_\tau(t)-\psi_\tau(t+\tau)}{\tau}u_\tau(t,x)\eta(x)\right.\right.\\&\quad+\left.\left.\psi_\tau(t)m\left(\left\lceil\frac{t}{\tau}\right\rceil\tau,u_\tau(t,x)\right)\left[\nabla u_\tau(t,x)f''(u_\tau(t,x))+\nabla\phi(x)\right]\cdot\nabla\eta(x)\right]\dd x \dd t\right|\le C\sqrt{\tau},
\end{split}
\end{align}
with the abbreviation $\psi_\tau(t):=\psi\left(\left\lceil\frac{t}{\tau}\right\rceil\tau\right)$.
\end{lemma}
\begin{proof}
Fix test functions $\eta\in C^\infty_c(\Rd)$ and $\psi\in C^\infty_c(\Rd)\cap C(\R_{\ge 0})$ and let $T>0$, $R>0$ such that $\supp\psi\subset [0,T]$ and $\supp\eta\subset \ball_R(0)$. We first consider $t=n\tau$ for a fixed $n\in\N$. In order to apply the flow interchange principle, we calculate the dissipation of $\ent$ along the semigroup $\flow_{(\cdot)}^{\aV_{n\tau}}(u_\tau^n)=:u_{(\cdot)}$:
\begin{align}
\label{eq:wdsp}
\begin{split}
-\frac{\dd}{\dd s}\ent(u_s)&=-\int_\Rd (f'(u_s)+\phi)(\sqrt{\tau}\Delta u_s+\dv(m(n\tau,u_s)\nabla\eta))\dd x\\
&=\sqrt{\tau}\int_\Rd(f''(u_s)|\nabla u_s|^2-u_s\Delta\phi)\dd x+\int_\Rd (f''(u_s)\nabla u_s+\nabla\phi)\cdot m(n\tau,u_s)\nabla\eta\dd x,
\end{split}
\end{align}
the last step obtained by integration by parts. We now are concerned with the passage to the limit as $s\searrow 0$. First, as $f''(z)\ge \gamma_0>0$ for all $z\ge 0$, one has
\begin{align*}
&\sqrt{\tau}\int_\Rd(f''(u_s)|\nabla u_s|^2-u_s\Delta\phi)\dd x+\int_\Rd (f''(u_s)\nabla u_s+\nabla\phi)\cdot m(n\tau,u_s)\nabla\eta\dd x\\
&=\sqrt{\tau}\int_\Rd(|g(u_s)|^2-u_s\Delta\phi)\dd x+\int_\Rd (\sqrt{f''(u_s)}g(u_s)+\nabla\phi)\cdot m(n\tau,u_s)\nabla\eta\dd x,
\end{align*}
defining the auxiliary map $g(u(x)):=\sqrt{f''(u(x))}\nabla u(x)$ for $u\in \X\cap H^1(\Rd)$ and $x\in\Rd$. Using Young's inequality, the estimates from Proposition \ref{prop:propm}(b) and $f''(z)\le \gamma_1$, we arrive at
\begin{align}\label{eq:wdspest}
\begin{split}
&-\frac{\dd}{\dd s}\ent(u_s)\\
&\ge \frac{\sqrt{\tau}}{2}\int_\Rd |g(u_s)|^2\dd x-\|\Delta \phi\|_{L^\infty}-C'(\|u_s\|_{L^2}^2+1)\\
&\ge \frac{\sqrt{\tau}}{2}\int_\Rd \gamma_0|\nabla u_s|^2\dd x-\|\Delta \phi\|_{L^\infty}-C'(\|u_s\|_{L^2}^2+1),
\end{split}
\end{align} 
for some constant $C'>0$. By similar arguments, one proves the boundedness of $\|u_s\|_{L^2}$ for small $s\ge 0$:
\begin{align}
\label{eq:L2bd}
\begin{split}
&-\frac{\dd }{\dd s}\int_\Rd u_s^2\dd x=2\sqrt{\tau}\int_\Rd|\nabla u_s|^2\dd x+2\int_\Rd m(n\tau,u_s)\nabla u_s\cdot\nabla\eta\dd x\\
&\ge \sqrt{\tau}\int_\Rd|\nabla u_s|^2\dd x-C'\int_\Rd u_s^2\dd x,
\end{split}
\end{align}
so Gronwall's lemma yields $\|u_s\|_{L^2}^2\le \|u_\tau^n\|_{L^2}^2\exp(C's)$. 

In view of the flow interchange lemma, we have by \eqref{eq:wdspest} that both $\|g(u_s)\|_{L^2}$ and $\|u_s\|_{H^1}$ are uniformly bounded for small $s\ge 0$. By Alaoglu's and Rellich's theorems, we infer (on a suitable subsequence)  that $u_s\rightharpoonup u_\tau^n$ in $H^1(\Rd)$ as well as $u_s\to u_\tau^n$ in $L^2(\ball_R(0))$ and $u_s(x)\to u_\tau^n(x)$ for a.e. $x\in\Rd$. Since $f''$ is uniformly bounded, the weak limit of $(g(u_s))_{s\ge 0}$ as $s\searrow 0$ in $L^2(\Rd)$ coincides with $g(u_\tau^n)$. Furthermore, since
\begin{align*}
\int_\Rd \left[m(n\tau,u_s)(\sqrt{f''(u_s)}+\nabla\phi)\cdot\nabla\eta\right]^2\dd x&\le C\int_{\ball_R(0)}(u_s^2+1)\dd x,
\end{align*}
and the integrand on the right-hand side is uniformly integrable thanks to $L^2$-convergence, one gets---using the weak convergence of $g(u_s)$ and the weak$\ast$-convergence of $u_s$:
\begin{align*}
&\lim_{s\searrow 0}\int_\Rd \left[-\sqrt{\tau}u_s\Delta\phi+\sqrt{f''(u_s)}g(u_s)+\nabla\phi)\cdot m(n\tau,u_s)\nabla\eta\right]\dd x\\
&=\int_\Rd \left[-\sqrt{\tau}u_\tau^n\Delta\phi+\sqrt{f''(u_\tau^n)}g(u_\tau^n)+\nabla\phi)\cdot m(n\tau,u_\tau^n)\nabla\eta\right]\dd x,
\end{align*}
by Vitali's convergence theorem. Combining this with the weak lower semicontinuity of the $L^2$ norm, we arrive at
\begin{align*}
&\dff^{\aV_{n\tau}}\ent(u_\tau^n)\ge\liminf_{s\searrow 0}\left(-\frac{\dd}{\dd s}\ent(u_s)\right)\\
&\ge \sqrt{\tau}\int_\Rd(f''(u_\tau^n)|\nabla u_\tau^n|^2-u_\tau^n\Delta\phi)\dd x+\int_\Rd (f''(u_\tau^n)\nabla u_\tau^n+\nabla\phi)\cdot m(n\tau,u_\tau^n)\nabla\eta\dd x.
\end{align*}
Now, the flow interchange lemma (Theorem \ref{thm:flowinterchange}) yields
\begin{align*}
&\aV_{n\tau}(u_\tau^n)+\tau\left[\sqrt{\tau}\int_\Rd(f''(u_\tau^n)|\nabla u_\tau^n|^2-u_\tau^n\Delta\phi)\dd x+\int_\Rd (f''(u_\tau^n)\nabla u_\tau^n+\nabla\phi)\cdot m(n\tau,u_\tau^n)\nabla\eta\dd x\right]\\
&-\frac{C}{2\sqrt{\tau}}\W_{m(n\tau,\cdot)}^2(u_\tau^n,u_\tau^{n-1})\le \aV_{n\tau}(u_\tau^{n-1}).
\end{align*}
Replacing $\eta$ with $-\eta$, we deduce the following chain of inequalities:
\begin{align}
&-\frac{C}{2\sqrt{\tau}}\W_{m(n\tau,\cdot)}^2(u_\tau^n,u_\tau^{n-1})+\sqrt{\tau}\left(\calH_{n\tau}(u_\tau^n)-\calH_{n\tau}(u_\tau^{n-1})\right)+\tau\sqrt{\tau}\int_\Rd(f''(u_\tau^n)|\nabla u_\tau^n|^2-u_\tau^n\Delta\phi)\dd x\nonumber\\
&\le \int_\Rd \eta(u_\tau^n-u_\tau^{n-1})\dd x+\tau\int_\Rd (f''(u_\tau^n)\nabla u_\tau^n+\nabla\phi)\cdot m(n\tau,u_\tau^n)\nabla\eta\dd x\label{eq:chainineq}\\
&\le \frac{C}{2\sqrt{\tau}}\W_{m(n\tau,\cdot)}^2(u_\tau^n,u_\tau^{n-1})-\sqrt{\tau}\left(\calH_{n\tau}(u_\tau^n)-\calH_{n\tau}(u_\tau^{n-1})\right)-\tau\sqrt{\tau}\int_\Rd(f''(u_\tau^n)|\nabla u_\tau^n|^2-u_\tau^n\Delta\phi)\dd x.\nonumber
\end{align}
Consider now the case of nonnegative-valued test functions $\psi$. Mutiply \eqref{eq:chainineq} with $\psi(n\tau)$ and sum up from $n=1$ to $n=\left\lfloor\frac{T}{\tau}\right\rfloor\tau+1=:N$ to arrive at
\begin{align}
\label{eq:chainineq2}
\begin{split}
&-\sum_{n=1}^N\psi(n\tau)\frac{C}{2\sqrt{\tau}}\W_{m(n\tau,\cdot)}^2(u_\tau^n,u_\tau^{n-1})+\sum_{n=1}^N\psi(n\tau)\sqrt{\tau}\left(\calH_{n\tau}(u_\tau^n)-\calH_{n\tau}(u_\tau^{n-1})\right)\\
&\qquad+\sum_{n=1}^N\psi(n\tau)\tau\sqrt{\tau}\int_\Rd(f''(u_\tau^n)|\nabla u_\tau^n|^2-u_\tau^n\Delta\phi)\dd x\\
&\le \tau\sum_{n=0}^{N-1}\frac{\psi(n\tau)-\psi((n+1)\tau)}{\tau}u_\tau^n\eta\dd x+\tau\sum_{n=1}^N\psi(n\tau)\int_\Rd (f''(u_\tau^n)\nabla u_\tau^n+\nabla\phi)\cdot m(n\tau,u_\tau^n)\nabla\eta\dd x\\
&\le \sum_{n=1}^N\psi(n\tau)\frac{C}{2\sqrt{\tau}}\W_{m(n\tau,\cdot)}^2(u_\tau^n,u_\tau^{n-1})-\sum_{n=1}^N\psi(n\tau)\sqrt{\tau}\left(\calH_{n\tau}(u_\tau^n)-\calH_{n\tau}(u_\tau^{n-1})\right)\\
&\qquad-\sum_{n=1}^N\psi(n\tau)\tau\sqrt{\tau}\int_\Rd(f''(u_\tau^n)|\nabla u_\tau^n|^2-u_\tau^n\Delta\phi)\dd x.
\end{split}
\end{align}
For a sign-changing test function $\psi$, we use \eqref{eq:chainineq2} for the positive and negative part, respectively, and subtract the resulting inequalities. Recalling that $|\psi|=\psi_++\psi_-$, we get
\begin{align}
\label{eq:chainineq3}
\begin{split}
&\left|\tau\sum_{n=0}^{N-1}\frac{\psi(n\tau)-\psi((n+1)\tau)}{\tau}u_\tau^n\eta\dd x+\tau\sum_{n=1}^N\psi(n\tau)\int_\Rd (f''(u_\tau^n)\nabla u_\tau^n+\nabla\phi)\cdot m(n\tau,u_\tau^n)\nabla\eta\dd x\right|\\
&\le \sum_{n=1}^N|\psi(n\tau)|\frac{C}{2\sqrt{\tau}}\W_{m(n\tau,\cdot)}^2(u_\tau^n,u_\tau^{n-1})+\left|\sum_{n=1}^N|\psi(n\tau)|\sqrt{\tau}\left(\calH_{n\tau}(u_\tau^n)-\calH_{n\tau}(u_\tau^{n-1})\right)\right|\\
&\qquad+\left|\sum_{n=1}^N|\psi(n\tau)|\tau\sqrt{\tau}\int_\Rd(f''(u_\tau^n)|\nabla u_\tau^n|^2-u_\tau^n\Delta\phi)\dd x\right|.
\end{split}
\end{align}
The left-hand side in \eqref{eq:chainineq3} coincides with the left-hand side of the desired estimate \eqref{eq:dweak1}---write in spatio-temporal integral form and recall the definition of $u_\tau$. It remains to control the right-hand side in \eqref{eq:chainineq3}.

First, thanks to \eqref{eq:distest} and Lemma \ref{lem:apri1}(c),
\begin{align*}
&\sum_{n=1}^N|\psi(n\tau)|\frac{C}{2\sqrt{\tau}}\W_{m(n\tau,\cdot)}^2(u_\tau^n,u_\tau^{n-1})+\left|\sum_{n=1}^N|\psi(n\tau)|\tau\sqrt{\tau}\int_\Rd(f''(u_\tau^n)|\nabla u_\tau^n|^2-u_\tau^n\Delta\phi)\dd x\right|\\
&\le \|\psi\|_{C^0}\sqrt{\tau}\left[C(\ent(u_0)-\inf_{u\in\X}\ent(u))+\gamma_1\|u_\tau\|_{L^2([0,T+\bar\tau];H^1)}^2+(T+\bar\tau)\|\Delta\phi\|_{L^\infty}\right]\le C'\sqrt{\tau},
\end{align*}
for a suitable constant $C'>0$. For the second part on the r.h.s. of \eqref{eq:chainineq3}, we first insert suitable terms and use the triangle inequality:
\begin{align*}
&\left|\sum_{n=1}^N|\psi(n\tau)|\left(\calH_{n\tau}(u_\tau^n)-\calH_{n\tau}(u_\tau^{n-1})\right)\right|\\
&\le \sum_{n=1}^N\big||\psi(n\tau)|-|\psi((n-1)\tau)|\big|\calH_{n\tau}(u_\tau^{n-1})+\sum_{n=1}^N|\psi((n-1)\tau)|\left|\calH_{n\tau}(u_\tau^{n-1})-\calH_{(n-1)\tau}(u_\tau^{n-1})\right|\\
&\quad+\left|\sum_{n=1}^N\left[|\psi((n-1)\tau)|\calH_{(n-1)\tau}(u_\tau^{n-1})-|\psi(n\tau)|\calH_{n\tau}(u_\tau^n)\right]\right|.
\end{align*}
From the regularity of $\psi$ and Lemma \ref{lem:heatent}(a), we have
\begin{align*}
&\sum_{n=1}^N\big||\psi(n\tau)|-|\psi((n-1)\tau)|\big|\calH_{n\tau}(u_\tau^{n-1})\\
&\le \sum_{n=1}^N\tau\|\psi\|_{C^1}C(1+\mom{u_\tau^n}+\|u_\tau^n\|_{L^2}^2)\\
&\le C(T+\bar\tau)\|\psi\|_{C^1}\sup_{t\in [0,T+\bar\tau]}\left(1+\mom{u_\tau(t)}+\|u_\tau\|_{L^\infty(L^2([0,T+\bar\tau];L^2))}^2\right),
\end{align*}
which is bounded uniformly in $\tau$ by Lemma \ref{lem:apri1}(a)\&(b). Second, the local Lipschitz condition for $\calH_t$ in the $t$ variable from (M4) yields
\begin{align*}
&\sum_{n=1}^N|\psi((n-1)\tau)|\left|\calH_{n\tau}(u_\tau^{n-1})-\calH_{(n-1)\tau}(u_\tau^{n-1})\right|\\
&\le \|\psi\|_{C^0}\tau\sum_{n=1}^N \int_\Rd\left(L_1(T+\bar\tau)(u_\tau^n)^{\beta_1(T+\bar\tau)}+L_2(T+\bar\tau)(u_\tau^n)^{\beta_2(T+\bar\tau)}\right)\dd x,
\end{align*}
which allows us to proceed as in the proof of \eqref{eq:carl} to obtain
\begin{align*}
&\sum_{n=1}^N|\psi((n-1)\tau)|\left|\calH_{n\tau}(u_\tau^{n-1})-\calH_{(n-1)\tau}(u_\tau^{n-1})\right|\le C.
\end{align*}
Finally, the last sum on the r.h.s. of \eqref{eq:chainineq3} is a telescopic sum and the remaining term can be controlled with \eqref{eq:carl} again:
\begin{align*}
&\left|\sum_{n=1}^N\left[|\psi((n-1)\tau)|\calH_{(n-1)\tau}(u_\tau^{n-1})-|\psi(n\tau)|\calH_{n\tau}(u_\tau^n)\right]\right|\\
&=|\psi(N\tau)\calH_{N\tau}(u_\tau^N)|\le C\|\psi\|_{C^0}\sup_{t\in [0,T+\bar\tau]}\left(1+\mom{u_\tau(t)}+\|u_\tau\|_{L^\infty(L^2([0,T+\bar\tau];L^2))}^2\right).
\end{align*}
All in all, we obtain the desired estimate \eqref{eq:dweak1}.
\end{proof}
By the same method, one considers the fourth order case (E2).
\begin{lemma}[Discrete weak formulation, case (E2)]\label{lem:dweak2}
Let $u_0\in \X(0)\cap H^1(\Rd)$, $\bar\tau>0$ and let test functions $\eta\in C^\infty_c(\Rd)$ and $\psi\in C^\infty_c((0,\infty))\cap C(\R_{\ge 0})$ be given. Then, there exists a constant $C>0$ such that for all $\tau\in (0,\bar\tau]$, one has
\begin{align}
\label{eq:dweak2}
\begin{split}
&\left|\int_0^\infty\int_\Rd \left[\frac{\psi_\tau(t)-\psi_\tau(t+\tau)}{\tau}u_\tau\eta+\psi_\tau q_\tau\cdot\nabla\phi\right.\right.\\&\quad+\left.\left.\psi_\tau
\trc\left\{\bigg(\Id_d\,\dv q_\tau \quad q_\tau\bigg)\nabla^2_{(p,z)}f(\nabla u_\tau,u_\tau)\begin{pmatrix}\nabla^2 u_\tau \\ (\nabla u_\tau)^\tT\end{pmatrix}\right\}\right]
\dd x \dd t\right|\le C\sqrt{\tau},
\end{split}
\end{align}
with the abbreviations $q_\tau:=m\left(\left\lceil\frac{t}{\tau}\right\rceil\tau,u_\tau\right)\nabla\eta$ and $\psi_\tau(t):=\psi\left(\left\lceil\frac{t}{\tau}\right\rceil\tau\right)$.
\end{lemma}
\begin{proof}
Again, we calculate the dissipation of $\ent$ along the semigroup $(\flow_s^{\aV_{n\tau}}(u_\tau^n))_{s\ge 0}$:
\begin{align*}
-\frac{\dd}{\dd s}\ent(u_s)&=-\int_\Rd\left[-\dv\nabla _pf(\nabla u_s,u_s)+\partial_z f(\nabla u_s,u_s)+\phi\right]\left[\sqrt{\tau}\Delta u_s+\dv(m(n\tau,u_s)\nabla\eta)\right]\dd x.
\end{align*}
Write $q_s:=m(n\tau,u_s)\nabla\eta$ for $s>0$ and $q_\tau^n:=m(n\tau,u_\tau^n)\nabla\eta$ for brevity. Integrating by parts as in the proof of Proposition \ref{prop:addreg1} and by elementary calculations, 
\begin{align*}
&-\frac{\dd}{\dd s}\ent(u_s)\\
&=\sqrt{\tau}\int_\Rd \left[\sum_{i=1}^d\begin{pmatrix}\partial_{x_i}\nabla u_s \\ \partial_{x_i}u_s\end{pmatrix}^\tT \nabla_{(p,z)}^2f(\nabla u_s,u_s)\begin{pmatrix}\partial_{x_i}\nabla u_s \\ \partial_{x_i}u_s\end{pmatrix}-u_s\Delta\phi\right]\dd x\\&\quad+\int_\Rd \left[\dv q_s\, \dv\nabla_p f(\nabla u_s,u_s)+q_s\cdot (\nabla\partial_z f(\nabla u_s,u_s)+\nabla\phi)\right]\dd x\\
&=\sqrt{\tau}\int_\Rd \left[\sum_{i=1}^d\begin{pmatrix}\partial_{x_i}\nabla u_s \\ \partial_{x_i}u_s\end{pmatrix}^\tT \nabla_{(p,z)}^2f(\nabla u_s,u_s)\begin{pmatrix}\partial_{x_i}\nabla u_s \\ \partial_{x_i}u_s\end{pmatrix}-u_s\Delta\phi\right]\dd x\\&\quad+\int_\Rd \left[\dv q_s\, \trc\left(\nabla^2 u_s\nabla^2_{p}f(\nabla u_s,u_s)+\nabla u_s(\nabla_p\partial_z f(\nabla u_s,u_s))^\tT\right)\right.\\&\quad\qquad+\left.q_s\cdot\left(\partial^2_z f(\nabla u_s,u_s)\nabla u_s+\nabla^2 u_s\nabla_p\partial_zf(\nabla u_s,u_s)+\nabla\phi\right)\right]\dd x\\
&=\sqrt{\tau}\int_\Rd \left[\sum_{i=1}^d\begin{pmatrix}\partial_{x_i}\nabla u_s \\ \partial_{x_i}u_s\end{pmatrix}^\tT \nabla_{(p,z)}^2f(\nabla u_s,u_s)\begin{pmatrix}\partial_{x_i}\nabla u_s \\ \partial_{x_i}u_s\end{pmatrix}-u_s\Delta\phi\right]\dd x+\int_\Rd q_s\cdot\nabla\phi\dd x\\&\quad+\int_\Rd \trc\left\{\bigg(\Id_d\,\dv q_s \quad q_s\bigg)\nabla^2_{(p,z)}f(\nabla u_s,u_s)\begin{pmatrix}\nabla^2 u_s \\ (\nabla u_s)^\tT\end{pmatrix}\right\}\dd x.
\end{align*}
Introducing the (uniquely defined) square root $(\nabla^2_{(p,z)}f(\nabla u(x),u(x)))^{1/2}$ of the positive definite matrix $\nabla^2_{(p,z)}f(\nabla u(x),u(x))$ and the vectors $g_i(u(x)):=(\nabla^2_{(p,z)}f(\nabla u(x),u(x)))^{1/2}\begin{pmatrix}\partial_{x_i}\nabla u(x) \\ \partial_{x_i}u(x)\end{pmatrix}$ for $i\in\{1,\ldots,d\}$, $x\in\Rd$ and $u\in\X\cap H^2(\Rd)$, the above can also be rewritten as
\begin{align*}
&-\frac{\dd}{\dd s}\ent(u_s)\\
&=\sqrt{\tau}\int_\Rd \left[\sum_{i=1}^d |g_i(u_s)|^2-u_s\Delta\phi\right]\dd x+\int_\Rd q_s\cdot\nabla\phi\dd x\\&\quad+\sum_{i=1}^d\int_\Rd \einsvec_i^\tT\bigg(\Id_d\,\dv q_s \quad q_s\bigg)(\nabla^2_{(p,z)}f(\nabla u_s,u_s))^{1/2}g_i(u_s)\dd x,
\end{align*}
where $\einsvec_i\in\Rd$ is the $i^\mathrm{th}$ canonical unit vector in $\Rd$. Using Young's inequality and the bounds from Proposition \ref{prop:propm}(b), we obtain
\begin{align}\label{eq:estq}
|q_s|&\le C(T+\bar\tau)\|\eta\|_{C^1}u_s,\\
|\dv q_s|&\le |\partial_z m(t,u_s)||\nabla u_s||\nabla\eta|+m(t,u_s)|\Delta\eta|\le C(T+\bar\tau)\|\eta\|_{C^2}(1+|u_s|+|\nabla u_s|),
\end{align}
and consequently
\begin{align*}
&-\frac{\dd}{\dd s}\ent(u_s)\\
&\ge \frac{\sqrt{\tau}}{2}\int_\Rd \sum_{i=1}^d |g_i(u_s)|^2\dd x-\sqrt{\tau}\|\Delta\phi\|_{L^\infty}-C(\|u_s\|_{L^2}^2+1)\\
&\quad-\frac12\sum_{i=1}^d\int_\Rd \left|\einsvec_i^\tT\bigg(\Id_d\,\dv q_s \quad q_s\bigg)(\nabla^2_{(p,z)}f(\nabla u_s,u_s))^{1/2}\right|^2\dd x\\
&\ge \frac{\sqrt{\tau}}{2}\int_\Rd \sum_{i=1}^d |g_i(u_s)|^2\dd x-C'(\|u_s\|_{H^1}^2+1)\\
&\ge \frac{\sqrt{\tau}}{2}\gamma_0 \|\nabla^2 u_s\|_{L^2}^2-C'(\|u_s\|_{H^1}^2+1),
\end{align*}
for a suitable constant $C'>0$. Furthermore, by similar arguments, one has
\begin{align*}
&-\frac{\dd }{\dd s}\int_\Rd |\nabla u_s|^2\dd x=\int_\Rd\left[\sqrt{\tau}(\Delta u_s)^2+\Delta u_s(\partial_z m(n\tau,u_s)\nabla u_s\cdot\nabla\eta+m(n\tau,u_s)\Delta\eta)\right]\dd x\\
&\ge \frac{\sqrt{\tau}}{2}\int_\Rd(\Delta u_s)^2\dd x-\frac12 \int_\Rd(\partial_z m(n\tau,u_s)\nabla u_s\cdot\nabla\eta+m(n\tau,u_s)\Delta\eta)^2\dd x\\
&\ge -C\left(\int_\Rd |\nabla u_s|^2\dd x+\|u_s\|_{L^2}+1\right)\\
&\ge -C\int_\Rd |\nabla u_s|^2\dd x-C',
\end{align*}
where we used in the last step that $\|u_s\|_{L^2}$ is uniformly bounded for small $s$, see \eqref{eq:L2bd}. Hence, Gronwall's lemma yields $s$-uniform boundedness of $u_s$ in $H^1(\Rd)$. All in all, we have shown that $u_s$ is bounded in $H^2$ as well as that the $g_i(u_s)$ are bounded in $L^2(\Rd)$. From now on, the argumentation for passing to the limit $s\searrow 0$ and deriving \eqref{eq:dweak2} is mutatis mutandis the same as in the proof of Lemma \ref{lem:dweak1}. We omit these technical details here.
\end{proof}
\subsection{Convergence}\label{subec:conv}
We are now in position to demonstrate the proofs of Theorems \ref{thm:exist2o}\&\ref{thm:exist4o}.\\

\noindent\textit{Proof of Theorem \ref{thm:exist2o}:} 
Let a vanishing sequence of step sizes $(\tau_k)_{k\in\N}$ be given and denote the associated sequence of discrete solutions by $(u_{\tau_k})_{k\in\N}$. The uniform estimate on the second moments and the $L^2([0,T];H^1(\Rd))$-norm from Lemma \ref{lem:apri1}(b)\&(c) yield in combination with the Alaoglu and Prokhorov theorems that $u_{\tau_k}(t)\rightharpoonup^\ast u(t)$ in $\prb$ as well as $u_{\tau_k}\rightharpoonup u$ in $L^2([0,T];H^1(\Rd))$, on a suitable non-relabelled subsequence as $k\to\infty$, for a limit map $u\in L^2([0,T];H^1(\Rd))$ with $u(t)\in\X$ for all $t\in[0,T]$. The additional H\"older continuity w.r.t. $\W_2$ is a consequence of the refined version of the Arzel\`{a}-Ascoli theorem \cite[Prop. 3.3.1]{savare2008} and \eqref{eq:holdest}. In order to prove strong convergence, we seek to apply Theorem \ref{thm:ex_aub}: choose $\ban=L^2(\ball_R(0))$ for fixed, but arbitrary $R>0$, the functional $\auxil:\,\ban\to[0,\infty]$ defined via
\begin{align*}
\auxil(u):=\begin{cases}\|u\|_{H^1}^2&\text{if }u\in H^1(\ball_R(0)),\\ +\infty&\text{otherwise},\end{cases}
\end{align*} 
and the pseudo-distance $\W:\ban\times\ban\to[0,\infty]$, defined by
\begin{align}\label{eq:pseudo}
\W(u,\tilde u):=\inf\left\{\right.\W_2(\rho,\tilde\rho)&:\,\rho,\tilde\rho\in\prb\text{ with}\\&\quad\mom{\rho}\le C_{T+\bar\tau},\,\mom{\tilde\rho}\le C_{T+\bar\tau},\,\rho|_{\ball_R(0)}=u,\,\tilde\rho|_{\ball_R(0)}=\tilde u\left.\right\},
\end{align}
where $C_{T+\bar\tau}$ is the constant from Lemma \ref{lem:apri1}(b) for $T+\bar\tau$ in place of $T$. Rellich's theorem immediately yields relative compactness of sublevels of $\auxil$. Concerning the pseudo-distance $\W$, the infimum (if it is finite) in its definition is always attained, since minimizing sequences belong to a weakly$\ast$-compact set by Prokhorov's theorem and the Wasserstein distance is lower semicontinuous w.r.t. weak$\ast$-convergence. One therefore easily shows that $\W(u,\tilde u)=0$ implies $u=\tilde u$ in $\ban$. For the lower semicontinuity property, let a sequence $(u_k,\tilde u_k)_{k\in\N}$ in $\ban\times\ban$ converging to $(u,\tilde u)\in\ban\times\ban$ be given. Without loss of generality, we may assume that $\W(u_k,\tilde u_k)<\infty$ for all $k\in\N$, hence there exists $(\rho_k,\tilde\rho_k)_{k\in\N}$ in $\prb\times\prb$ with $\W(u_k,\tilde u_k)=\W_2(\rho_k,\tilde\rho_k)$ for all $k$. Clearly, on a non-relabelled subsequence, $\rho_k\rightharpoonup^\ast \rho'$ and $\tilde\rho_k\rightharpoonup^\ast \tilde\rho'$ in $\prb$, so
\begin{align*}
\W(u,\tilde u)\le \W_2(\rho',\tilde\rho')\le \liminf_{k\to\infty}\W_2(\rho_k,\tilde\rho_k)=\liminf_{k\to\infty}\W(u_k,\tilde u_k),
\end{align*}
which was to be proved. Henceforth, one has for all $t\ge 0$, $h>0$ and $k\in\N$ that
\begin{align}\label{eq:pseudoab}
\W(u_{\tau_k}|_{\ball_R(0)}(t+h),u_{\tau_k}|_{\ball_R(0)}(t))&\le \W_2(u_{\tau_k}(t+h),u_{\tau_k}(t)).
\end{align}
Hypothesis \eqref{eq:hypo1} from Theorem \ref{thm:ex_aub} is obviously satisfied thanks to Lemma \ref{lem:apri1}(c). Concerning hypothesis \eqref{eq:hypo2} from Theorem \ref{thm:ex_aub}, we claim that
\begin{align}
\label{eq:hypo2ver}
\begin{split}
&\sup_{k\in\N}\int_0^{T-h}\W_2(u_{\tau_k}(t+h),u_{\tau_k}(t))\dd t\\&\le \max\left(1,\sqrt{T+\bar\tau}\right)\sqrt{2C(T+\bar\tau)(\ent(u_0)-\inf\ent)(T+\bar \tau)h}, 
\end{split}
\end{align}
for all $h\in(0,\bar\tau)$, from which \eqref{eq:hypo2} follows using \eqref{eq:pseudoab}. Indeed, for fixed $k\in\N$ and $h\in (0,\tau_k]$, one has, for $N:=\gau{\frac{T}{\tau_k}}$,
\begin{align*}
&\int_0^{T-h}\W_2(u_{\tau_k}(t+h),u_{\tau_k}(t))\dd t=\sqrt{C(T+\bar\tau)}\sum_{n=1}^{N}h\W_{m(n\tau_k,\cdot)}(u_{\tau_k}^n,u_{\tau_k}^{n+1})\\&\le \sqrt{2C(T+\bar\tau)(\ent(u_0)-\inf\ent)}\sqrt{h^2 N}\\
&\le \sqrt{2(\ent(u_0)-\inf\ent)(T+\bar\tau)h},
\end{align*}
thanks to H\"older's inequality and the total square distance estimate \eqref{eq:distest}. On the other hand, for $h\in(\tau_k,\bar\tau]$, we directly get from the H\"older type estimate \eqref{eq:holdest}:
\begin{align*}
\int_0^{T-h}\W_2(u_{\tau_k}(t+h),u_{\tau_k}(t))\dd t&\le (T-h)\sqrt{2C(T+\bar\tau)(\ent(u_0)-\inf\ent)h}\\&\le (T+\bar\tau)\sqrt{2C(T+\bar\tau)(\ent(u_0)-\inf\ent)h}.
\end{align*}
Theorem \ref{thm:ex_aub} thus is applicable for the family $(u_{\tau_k}|_{\ball_R(0)})_{k\in\N}$; hence, $u_{\tau_k}(t)$ converges to $u(t)$ in $L^2(\ball_R(0))$, in measure w.r.t. $t\in [0,T]$. By the dominated convergence theorem and the uniform estimate in $L^\infty([0,T];L^2(\Rd))$ from Lemma \ref{lem:apri1}(a), $u_{\tau_k}$ converges to $u$ in $L^2([0,T]\times\ball_R(0))$, and $u\in L^\infty([0,T];L^2(\Rd))$. After a diagonal argument (considering a family of radii $R\uparrow +\infty$), extracting further subsequences, we have that 
\begin{align}
u_{\tau_k}\to u&\quad\text{in }L^2([0,T];L^2(\ball_R(0))\text{ for all $R>0$ simultaneously},\label{eq:L2tx}\\
u_{\tau_k}(t)\to u(t)&\quad\text{in }L^2(\ball_R(0))\text{ for all $R>0$ simultaneously and almost all $t\in [0,T]$},\label{eq:L2x}\\
u_{\tau_k}(t,x)\to u(t,x)&\quad\text{for almost all }(t,x)\in [0,T]\times\Rd.\label{eq:aec}
\end{align}
Fix now a $t\in[0,T]$ such that \eqref{eq:L2x} holds. Since by construction of $u_{\tau_k}$, one has
\begin{align*}
u_{\tau_k}(t,x)\le S\left(\left\lceil\frac{t}{\tau_k}\right\rceil\tau_k\right),
\end{align*}
for almost all $x\in\Rd$, and $S$ is continuous thanks to condition (M2), one also has
\begin{align*}
u(t,x)-S(t)\le \lim_{k\to\infty}\left(u(t,x)-u_{\tau_k}(t,x)+S\left(\left\lceil\frac{t}{\tau_k}\right\rceil\tau_k\right)-S(t)\right)=0,
\end{align*}
for almost all $x\in\Rd$. Hence, $u(t,x)\in[0,S(t)]$ a.e.

It now remains to prove that the limit map $u$ is a solution to \eqref{eq:pde2o} in the sense of distributions, i.e.: For all $\eta\in C^\infty_c(\Rd)$ and all $\psi\in C^\infty_c((0,\infty))\cap C(\R_{\ge 0})$, it holds that
\begin{align}
\label{eq:cweak1}
\begin{split}
&\int_0^\infty\int_\Rd \left[-\partial_t\psi(t) u(t,x)\eta(x)+\psi(t)m(t,u(t,x))\left[\nabla u(t,x) f''(u(t,x))+\nabla\phi(x)\right]\cdot\nabla\eta(x)\right]\dd x \dd t\\&\qquad=0,
\end{split}
\end{align}
Since $t\mapsto\frac{\psi_{\tau_k}(t)-\psi_{\tau_k}(t+{\tau_k})}{\tau_k}$ converges uniformly in $\R_{\ge 0}$ to $-\partial_t\psi$ and $u_{{\tau_k}}$ converges weakly to $u$, we immediately have that 
\begin{align*}
\lim_{k\to\infty}\int_0^\infty\int_\Rd \frac{\psi_{\tau_k}(t)-\psi_{\tau_k}(t+{\tau_k})}{\tau_k} u_{\tau_k}(t,x)\eta(x)\dd x\dd t&=\int_0^\infty\int_\Rd -\partial_t\psi(t) u(t,x)\eta(x)\dd x\dd t.
\end{align*}
Thus, by the discrete weak formulation \eqref{eq:dweak1}, it remains to show that (possibly extracting another subsequence) one has
\begin{align}\label{eq:nonlin}\begin{split}
\lim_{k\to\infty}\int_0^\infty\int_\Rd \bigg[&\psi_{\tau_k}(t)m\left(\left\lceil\frac{t}{{\tau_k}}\right\rceil{\tau_k},u_{\tau_k}(t,x)\right)\left[\nabla u_{\tau_k}(t,x)f''(u_{\tau_k}(t,x))+\nabla\phi(x)\right]\cdot\nabla\eta(x)\bigg.\\
&-\psi(t)m(t,u(t,x))\left[\nabla u(t,x) f''(u(t,x))+\nabla\phi(x)\right]\cdot\nabla\eta(x)\bigg.\bigg]\dd x \dd t=0.
\end{split}\end{align}
First, thanks to the convergence almost everywhere from \eqref{eq:aec} and the continuity of $m$ and $f''$, 
\begin{align*}
&\psi_{\tau_k}(t)m\left(\left\lceil\frac{t}{{\tau_k}}\right\rceil{\tau_k},u_{\tau_k}(t,x)\right)\left[f''(u_{\tau_k}(t,x))+\nabla\phi(x)\right]\cdot\nabla\eta(x)\\&\quad-\psi(t)m(t,u(t,x))\left[f''(u(t,x))+\nabla\phi(x)\right]\cdot\nabla\eta(x)\stackrel{k\to\infty}{\longrightarrow} 0,
\end{align*}
a.e. on $[0,T]\times \ball_R(0)\supset \supp\psi\times\supp\eta$. Furthermore, due to the bounds on $f''$ and $m$ (recall Proposition \ref{prop:propm}(b)), one has
\begin{align*}
&\int_0^\infty\int_\Rd \bigg|\bigg.\psi_{\tau_k}(t)m\left(\left\lceil\frac{t}{{\tau_k}}\right\rceil{\tau_k},u_{\tau_k}(t,x)\right)\left[f''(u_{\tau_k}(t,x))+\nabla\phi(x)\right]\cdot\nabla\eta(x)\\
&\qquad-\psi(t)m(t,u(t,x))\left[f''(u(t,x))+\nabla\phi(x)\right]\cdot\nabla\eta(x)\bigg.\bigg|^2\\
&\le C\int_0^\infty\int_\Rd \eins{\supp\psi\times\supp\eta}(t,x)(1+u_{\tau_k}(t,x)^2+u(t,x)^2)\dd x\dd t,
\end{align*}
for some $C>0$, and the right-hand side is uniformly integrable since $u_{\tau_k}\to u$ in $L^2([0,T];L^2(\ball_\R(0)))$, cf. \eqref{eq:L2tx}. Thus, Vitali's convergence theorem yields that
\begin{align*}
&\int_0^\infty\int_\Rd \bigg|\bigg.\psi_{\tau_k}(t)m\left(\left\lceil\frac{t}{{\tau_k}}\right\rceil{\tau_k},u_{\tau_k}(t,x)\right)\left[f''(u_{\tau_k}(t,x))+\nabla\phi(x)\right]\cdot\nabla\eta(x)\\
&\qquad-\psi(t)m(t,u(t,x))\left[f''(u(t,x))+\nabla\phi(x)\right]\cdot\nabla\eta(x)\bigg.\bigg|^2\stackrel{k\to\infty}{\longrightarrow}0.
\end{align*}
Recalling the weak convergence of $\nabla u_{\tau_k}$ to $\nabla u$ in $L^2([0,T]\times\Rd)$, we consequently infer \eqref{eq:nonlin} and also the desired distributional formulation \eqref{eq:cweak1}.
\hfill\qed\\

The proof of Theorem \ref{thm:exist4o} goes along similar lines; we now highlight the important differences.\\

\noindent\textit{Proof of Theorem \ref{thm:exist4o}:} 
Our strategy for proving strong convergence of $u_{\tau_k}\to u$ in the space \break $L^2([0,T];H^1(\ball_R(0)))$ is adapted by choosing the space $\ban:=H^1(\ball_R(0))$ and the functional $\auxil:\,\ban\to[0,\infty]$ defined via
\begin{align*}
\auxil(u):=\begin{cases}\|u\|_{H^2}^2&\text{if }u\in H^2(\ball_R(0)),\\ +\infty&\text{otherwise}.\end{cases}
\end{align*} 
The main arguments stay---mutatis mutandis---the same. To show that $u$ is a weak solution to \eqref{eq:pde4o}, one needs to verify that
\begin{align*}
\int_0^\infty\int_\Rd\bigg[\bigg.&\psi_{\tau_k}\trc\left\{\bigg(\Id_d\,\dv q_{\tau_k} \quad q_{\tau_k}\bigg)\nabla^2_{(p,z)}f(\nabla u_{\tau_k},u_{\tau_k})\begin{pmatrix}\nabla^2 u_{\tau_k} \\ (\nabla u_{\tau_k})^\tT\end{pmatrix}\right\}\\
&-\psi\trc\left\{\bigg(\Id_d\,\dv q \quad q\bigg)\nabla^2_{(p,z)}f(\nabla u,u)\begin{pmatrix}\nabla^2 u \\ (\nabla u)^\tT\end{pmatrix}\right\}\bigg.\bigg]\dd x \dd t\stackrel{k\to\infty}{\longrightarrow}0,
\end{align*}
for $q:=m(t,u)\nabla\eta$, or that, equivalently, for each $i\in\{1,\ldots,d\}$:
\begin{align*}
\int_0^\infty\int_\Rd\bigg[\bigg.&\psi_{\tau_k}\einsvec_i\bigg(\Id_d\,\dv q_{\tau_k} \quad q_{\tau_k}\bigg)\nabla^2_{(p,z)}f(\nabla u_{\tau_k},u_{\tau_k})\begin{pmatrix}\partial_{x_i}\nabla u_{\tau_k} \\ \partial_{x_i} u_{\tau_k}\end{pmatrix}\\
&-\psi\einsvec_i\bigg(\Id_d\,\dv q \quad q\bigg)\nabla^2_{(p,z)}f(\nabla u,u)\begin{pmatrix}\partial_{x_i}\nabla u \\ \partial_{x_i}\nabla u\end{pmatrix}\bigg.\bigg]\dd x \dd t\stackrel{k\to\infty}{\longrightarrow}0,
\end{align*}
In view of the weak convergence of $u_{\tau_k}\rightharpoonup u$ in $L^2([0,T];H^2(\Rd))$, it suffices to show that
\begin{align}\label{eq:nonlin4o}
\psi_{\tau_k}\einsvec_i\bigg(\Id_d\,\dv q_{\tau_k} \quad q_{\tau_k}\bigg)\nabla^2_{(p,z)}f(\nabla u_{\tau_k},u_{\tau_k})\stackrel{k\to\infty}{\longrightarrow}\psi\einsvec_i\bigg(\Id_d\,\dv q \quad q\bigg)\nabla^2_{(p,z)}f(\nabla u,u)
\end{align}
in $L^2([0,T];L^2(\Rd;\Rd))$. Thanks to the continuity of $f$, $m$ and $\partial_z m$ and the strong convergence of $u_{\tau_k}\to u$ in $L^2([0,T];H^2(\ball_R(0)))$ for $\ball_R(0)\supset\supp\eta$, \eqref{eq:nonlin4o} holds pointwise almost everywhere on $[0,T]\times\ball_R(0)$, since 
\begin{align*}
\dv q_{\tau_k}&=\partial_z m\left(\left\lceil\frac{t}{\tau}\right\rceil\tau_k,u_{\tau_k}\right)\nabla u_{\tau_k}\nabla\eta+m\left(\left\lceil\frac{t}{\tau}\right\rceil\tau_k,u_{\tau_k}\right)\Delta\eta.
\end{align*} 
By employing the elementary bounds on $m$ and $\partial_z m$ from Proposition \ref{prop:propm}(b), we deduce bounds for $q_{\tau_k}$ analogous to those in \eqref{eq:estq}. Consequently, recalling the bound on $\nabla_{(p,z)}^2 f$ from (E2), we get
\begin{align*}
&\int_0^\infty\int_\Rd \left|\psi_{\tau_k}\einsvec_i\bigg(\Id_d\,\dv q_{\tau_k} \quad q_{\tau_k}\bigg)\nabla^2_{(p,z)}f(\nabla u_{\tau_k},u_{\tau_k})\right|^2\dd x\dd t\\
&\le C\int_0^\infty\int_\Rd \eins{\supp\psi\times\supp\eta}(t,x)(1+u_{\tau_k}(t,x)^2+|\nabla u_{\tau_k}(t,x)|^2)\dd x\dd t,
\end{align*}
and the right-hand side is uniformly integrable. Vitali's theorem thus yields the desired convergence \eqref{eq:nonlin4o} in $L^2([0,T];L^2(\Rd;\Rd))$.
\hfill\qed\\
\section{Generalizations}\label{sec:gen}
We give some concluding remarks about possible generalizations of our Theorems \ref{thm:exist2o} and \ref{thm:exist4o}.
\subsection{Bounded spatial domains}\label{subsec:bdd}
In principle, our results also hold if the spatial variable $x$ is confined to a bounded and convex domain $\Omega\subset\R^d$ with sufficiently smooth boundary $\partial\Omega$, when we consider the initial-boundary-value problem for \eqref{eq:pde} with the no-flux and Neumann boundary conditions
\begin{align}\label{eq:bc}
m(t,u(t,\bar x))\partial_\nu \frac{\delta\ent}{\delta u}(u(t,\bar x))&=0=\partial_\nu u(t,\bar x),
\end{align}
for all $t>0$ and $\bar x\in\partial\Omega$, where $\nu$ denotes the unit normal vector field to $\partial\Omega$.

In this framework, one may replace the conditions (E1) and (E2) on the free energy functional by the following different, mostly weaker assumptions.
\begin{definition}[Admissible energy functionals on bounded spatial domains]~
\begin{enumerate}[{(E}1{')}]
\item \emph{Second order equations \eqref{eq:pde2o}}: Let $\ent$ be of the form
\begin{align*}
\ent(u)=\int_\Omega \left[f(u(x))+\phi(x)u(x)\right]\dd x,
\end{align*}
where $\phi\in C^2(\bar\Omega)$ and $f\in C^2(\R_{\ge 0})$ with $f''(z)\in [\gamma_0,\gamma_1]$ for $0<\gamma_0\le \gamma_1<\infty$ and all $z\ge 0$.
\item \emph{Fourth order equations \eqref{eq:pde4o}}: Let $\ent$ be of the form
\begin{align*}
\ent(u)=\int_\Omega \left[f(\nabla u(x),u(x))+\phi(x)u(x)\right]\dd x,
\end{align*}
with $\phi$ as in (E1'). For $f\in C^2(\Rd\times \R_{\ge 0})$, we assume the following:

There exist $0<\gamma_0\le \gamma_1<\infty$ such that for each $(p,z)\in\Rd\times\R_{\ge 0}$ and all $(v_p,v_z)\in \R^d\times\R$, one has
\begin{align*}
\gamma_0|v_p|^2&\le\begin{pmatrix}v_p\\ v_z\end{pmatrix}^\tT\nabla^2_{(p,z)}f(p,z)\begin{pmatrix}v_p\\ v_z\end{pmatrix}\le \gamma_1(|v_p|^2+v_z^2),
\end{align*}
i.e. $f$ is convex w.r.t. $(p,z)$ and uniformly convex w.r.t. $p$. 

Additionally, we require the following \emph{isotropy condition} w.r.t. the gradient variable $p$: There exists a map $\tilde f\in C^2(\R_{\ge 0}\times\R_{\ge 0})$ with $\partial_r \tilde f(r,z)\ge 0$ for all $(r,z)\in \R_{\ge 0}\times\R_{\ge 0}$ such that $f(p,z)=\tilde f(|p|,z)$ for all $(p,z)\in\Rd\times \R_{\ge 0}$.
\end{enumerate}
\end{definition}
Notice that one can in the case (E2') now also allow for integrands $f$ which do not depend on $z$. Since $\Omega$ is bounded, a control on the gradient $\nabla u$ in $L^2$ also yields bounds in $L^2$ on $u$ itself via Poincar\'{e}'s inequality. Clearly, we also do not have to require finiteness of second moments $\mom{u}$ for probability measures $u\in\prb(\Omega)$ on $\Omega$ explicitly anymore. However, due the appearance of boundary terms, we need the additional isotropy condition (compare to \cite{loibl2015}).

The respective existence results---which can be obtained by essentially the same methods as Theorems \ref{thm:exist2o} and \ref{thm:exist4o}---read as follows.
\begin{thm}[Existence: second order case (E1')]\label{thm:exist2obd}
Assume that the mobility $m$ satisfies \textup{(M1)--(M4)} and the energy functional $\ent$ is of the form \textup{(E1')} and let an initial datum $u_0\in \X\cap L^2(\Rd)$ with $u_0(x)\in [0,S(0)]$ for almost every $x\in\Omega$ be given. Then, for each $\tau>0$, the map $u_\tau$ obtained via \eqref{eq:mms} is well-defined. Furthermore, for each vanishing sequence $\tau_k\searrow 0$ $(k\to\infty)$, there exists a (nonrelabelled) subsequence and a limit map $u:\R_{\ge 0}\to\X$ such that the following is true for each fixed $T>0$:
\begin{enumerate}[(a)]
\item $u\in L^\infty([0,T];L^2(\Omega))\cap L^2([0,T];H^1(\Omega))\cap C^{1/2}([0,T];\W_2)$, and for almost every $t\ge 0$, one has $u(t,x)\le S(t)$ for almost all $x\in\Omega$;
\item $u_{\tau_k}(t)\to u(t)$ uniformly in $t\in[0,T]$ with respect to the distance $\W_2$ on $\prb(\Omega)$;
\item $u_{\tau_k}$ converges to $u$ strongly in $L^2([0,T];L^2(\Omega))$ and weakly in $L^2([0,T];H^1(\Omega))$;
\item $u$ is a solution to \eqref{eq:pde2o} in the sense of distributions, it attains the boundary conditions \eqref{eq:bc} and the initial condition $u(0,\cdot)=u_0$ holds a.e. on $\Omega$.
\end{enumerate}
\end{thm}
In the fourth order case, we have:
\begin{thm}[Existence: fourth order case (E2')]\label{thm:exist4obd}
Assume that the mobility $m$ satisfies \textup{(M1)--(M4)} and the energy functional $\ent$ is of the form \textup{(E2')} and let an initial datum $u_0\in \X\cap H^1(\Omega)$ with $u_0(x)\in [0,S(0)]$ for almost every $x\in\Omega$ be given. Then, for each $\tau>0$, the map $u_\tau$ obtained via \eqref{eq:mms} is well-defined. Furthermore, for each vanishing sequence $\tau_k\searrow 0$ $(k\to\infty)$, there exists a (nonrelabelled) subsequence and a limit map $u:\R_{\ge 0}\to\X$ such that the following is true for each fixed $T>0$:
\begin{enumerate}[(a)]
\item $u\in L^\infty([0,T];H^1(\Omega))\cap L^2([0,T];H^2(\Omega))\cap C^{1/2}([0,T];\W_2)$, and for almost every $t\ge 0$, one has $u(t,x)\le S(t)$ for almost all $x\in\Omega$;
\item $u_{\tau_k}(t)\to u(t)$ uniformly in $t\in[0,T]$ with respect to the distance $\W_2$ on $\prb(\Omega)$;
\item $u_{\tau_k}$ converges to $u$ strongly in $L^2([0,T];H^1(\Omega))$ and weakly in $L^2([0,T];H^2(\Omega))$;
\item $u$ is a solution to \eqref{eq:pde2o} in the sense of distributions, it attains the boundary conditions \eqref{eq:bc} and the initial condition $u(0,\cdot)=u_0$ holds a.e. on $\Omega$.
\end{enumerate}
\end{thm}
\subsection{Non-Lipschitz mobility functions}\label{subsec:nonlip}
A second possibility of extension is concerned with the assumptions (M1)--(M4) on the mobility function $m$. Specifically, we consider mobilities which satisfy (M1),(M2) and (M4) but not the Lipschitz-semiconcavity condition (M3). In this framework, we restrict ourselves to the case of a bounded spatial domain $\Omega\subset\Rd$, as in Section \ref{subsec:bdd}. Our strategy is as follows: We approximate $m$ similarly to \cite{lisini2012} with mobilities $m_\delta$---which satisfy (M3)---for small $\delta>0$ (see Definition \ref{def:appr} below). The respective family of weak solutions $(u_\delta)_{\delta>0}$ provided by the Theorems \ref{thm:exist2obd} and \ref{thm:exist4obd} for $m_\delta$ in place of $m$ is then expected to converge to a limit $u$ which is a weak solution of the original problem. The precise assumptions read as follows:
\begin{definition}[Non-Lipschitz mobilities]~
\begin{enumerate}[{(M}1{')}]
\item $m$ satisfies condition (M1) for a \emph{constant} function $S(\cdot)$ (with value $S\in \R_{>0}\cup\{+\infty\}$).
\item Condition (M2) holds.
\item There exists a map $\tilde m: [0,S]\to \R$ satisfying (M1) with the property
\begin{align*}
m(t,z)&\le \tilde m(z)\quad\text{for all $t\ge 0$ and all $z\in [0,S]$.}
\end{align*}
\item Condition (M4) holds.
\end{enumerate}
When considering fourth-order equations (cf.~(E2')), due to the appearance of gradients, the following condition on $\partial_z m$ is needed:
\begin{enumerate}[{(M}1{')}]
\setcounter{enumi}{4}
\item One has for each $t\ge 0$ that
\begin{align*}
\lim_{z\searrow 0}\partial_z m(t,z)\sqrt{z}=0,
\end{align*}
and, if $S<\infty$, also
\begin{align*}
\lim_{z\nearrow S}\partial_z m(t,z)\sqrt{S-z}=0,
\end{align*}
\end{enumerate}
\end{definition}
Notice that we do not require that the auxiliary time-independent function $\tilde m$ satisfies (M3): The distance $\W_{\tilde m}$ acts as a surroggate for the $L^2$-Wasserstein distance $\W_2$ for the H\"older type estimate \eqref{eq:holdest} which is needed in the proof of strong convergence. We restricted our considerations on constant value spaces $[0,S]$ mostly due to technical reasons: for the approximation in Definition \ref{def:appr} below, it is convenient to have equal support for all $m(t,\cdot)$. Condition (M5') in particular ensures the convergence of $\nabla m(t,u_\delta)$ to the respective limit; using the well-known \emph{Lions-Villani estimate} on square roots (see \cite{lisini2012} for the detailed calculation) provides the necessary a priori estimate.\\

The following examples are paradigmatic; they correspond to time-dependent versions of the examples discussed in \cite{dns2009,lisini2010,lisini2012}.
\begin{example}[Paradigmatic choices]~\label{ex:para2}
\begin{enumerate}[(a)]
\item For $S=+\infty$ and a sufficiently regular map $\alpha$ with values in $(\frac12,1]$, set
\begin{align*}
m(t,z)&=z^{\alpha(t)}\qquad\text{for all $t\ge 0$ and $z\ge 0$.}
\end{align*}
Choosing the \emph{Dirichlet} energy
\begin{align}\label{eq:diri}
\ent_D(u)=\frac12\int_\Omega |\nabla u(x)|^2\dd x,
\end{align}
which now is admissible in (E2'), \eqref{eq:pde4o} reads as the thin film equation with time-dependent mobility:
\begin{align*}
\partial_t u&=-\dv(u^{\alpha(t)}\nabla\Delta u).
\end{align*}
\item For $S<+\infty$ and sufficiently regular maps $\sigma_1$, $\sigma_2$ with values in $(\frac12,1]$, set
\begin{align*}
m(t,z)&=z^{\sigma_1(t)}(S-z^{\sigma_2(t)})\qquad\text{for all $t\ge 0$ and $z\in [0,S]$.}
\end{align*}
\end{enumerate}
\end{example}
With the conditions (M1')--(M4'), we can---uniformly with respect to time $t$---define the approximating mobilities $m_\delta$ as in \cite{lisini2012} distinguishing the cases $S<\infty$ and $S=+\infty$:
\begin{definition}[Approximation of $m$]\label{def:appr}
Let $\delta\in (0,\bar\delta)$ and $\bar\delta>0$ sufficiently small. \\
If $S<\infty$, define for each fixed $t\ge 0$ and $z\in [0,S]$:
\begin{align}\label{eq:app1}
m_\delta(t,z):=m\left(t,\frac{z_{\delta,1}(t)-z_{\delta,2}(t)}{S}z-z_{\delta,1}(t)\right)-\delta,
\end{align}
where $z_{\delta,1}(t)<z_{\delta,2}(t)$ are the two solutions of $m(t,z)=\delta$ for $z$ at fixed $t$.\\
If $S=\infty$, define
\begin{align}\label{eq:app2}
m_\delta(t,z):=m(t,z+z_\delta(t))-\delta,
\end{align}
for all $t\ge 0$ and $z\ge 0$,
where $z_\delta(t)$ is the unique solution of $m(t,z)=\delta$ for $z$ at fixed $t$.
\end{definition}
Thanks to the assumed regularity of $m$ w.r.t. $t$ and the fact that $S(\cdot)$ remains constant, elementary properties of the approximation in Definiton \ref{def:appr} carry over from the autonomous case \cite{lisini2012} to the non-autonomous case. One has (compare to \cite{lisini2012}):
\begin{lemma}[Properties of the approximation]\label{lem:propappr}
Let $\bar\delta>0$ be sufficiently small. The following statements hold:
\begin{enumerate}[(a)]
\item For each $\delta\in(0,\bar\delta)$, the map $m_\delta$ satisfies \textup{(M1)--(M4)} and the estimate $0\le m_\delta(t,z)\le m(t,z)$ for all $t\ge 0$ and $z\in [0,S]$.
\item For each $\delta\in(0,\bar\delta)$, all $t\ge 0$ and all $u,\tilde u\in \X(0)$:
\begin{align*}
\W_{\tilde m}(u,\tilde u)\le \W_{m(t,\cdot)}(u,\tilde u)\le \W_{m_\delta(t,\cdot)}(u,\tilde u).
\end{align*}
\item As $\delta\searrow 0$, $m_\delta\to m$ locally uniformly in $\R_{\ge 0}\times [0,S]$ and $\partial_z m_\delta\to \partial_z m$ locally uniformly in $\R_{\ge 0}\times (0,S)$.
\end{enumerate}
\end{lemma}
Using the methods from \cite{lisini2012} and \cite{zinsl2016}, one finds thanks to Lemma \ref{lem:propappr} that the map $h_\delta$ from assumption (M4) admits for all $T>0$ the bound
\begin{align*}
|h_\delta(t,z)|\le C(z^2+1)\quad\text{for all $t\in[0,T]$ and $z\in [0,S]$},
\end{align*}
for a constant $C>0$ which is \emph{independent} from $\delta\in (0,\bar\delta)$. In consequence, the following $\delta$-uniform \emph{a priori} estimates hold:
\begin{lemma}[$\delta$-uniform a priori estimates, case (E1')]
Let $\bar\delta>0$ sufficiently small and let $T>0$. Then, there exists a constant $C>0$ such that for all $\delta\in(0,\bar\delta)$ and all $t\in[0,T]$:
\begin{align*}
\|u_\delta\|_{L^\infty([0,T];L^2(\Rd))}&\le C,\\
\|u_\delta\|_{L^2([0,T];H^1(\Rd))}&\le C,\\
\W_{\tilde m}(u_\delta(t),u_\delta(s))&\le C\sqrt{|t-s|}\quad\text{for all }s,t\in[0,T].
\end{align*}
\end{lemma}
\begin{lemma}[$\delta$-uniform a priori estimates, case (E2')]
Let $\bar\delta>0$ sufficiently small and let $T>0$. Then, there exists a constant $C>0$ such that for all $\delta\in(0,\bar\delta)$ and all $t\in[0,T]$:
\begin{align*}
\|u_\delta\|_{L^\infty([0,T];H^1(\Rd))}&\le C,\\
\|u_\delta\|_{L^2([0,T];H^2(\Rd))}&\le C,\\
\W_{\tilde m}(u_\delta(t),u_\delta(s))&\le C\sqrt{|t-s|}\quad\text{for all }s,t\in[0,T].
\end{align*}
\end{lemma}
With similar arguments as in Section \ref{subec:conv}, one passes to the limit $u_\delta\to u$ as $\delta\searrow 0$. The required strong convergence can again be achieved via Theorem \ref{thm:ex_aub}: the necessary pseudo-distance $\W$ is chosen to be the (auxiliary) distance $\W_{\tilde m}$. Note that since we work in bounded spatial domains $\Omega\subset\Rd$, an infimization process as in \eqref{eq:pseudo} is not needed here. As a last step in the proof, one passes to the limit in the weak formulation of \eqref{eq:pde2o} or \eqref{eq:pde4o} for $m_\delta$ and $u_\delta$ in place of $m$ and $u$, respectively. Again, similar methods as in Section \ref{subec:conv} apply; convergence of the mobilities and the mobility gradients can be proved by essentially the same strategy as in \cite{lisini2012}. In the end, one arrives at:
\begin{thm}[Existence for non-Lipschitz mobilities: second order case (E1')]\label{thm:exist2onl}
Assume that the \break mobility $m$ satisfies \textup{(M1')--(M4')} and that the energy functional $\ent$ is of the form \textup{(E1')}. Let an initial datum $u_0\in \X\cap L^2(\Rd)$ with $u_0(x)\in [0,S]$ for almost every $x\in\Omega$ be given. Define $m_\delta$ for $\delta\in (0,\bar\delta)$ and sufficiently small $\bar\delta$ as in Definition \ref{def:appr} and let $u_\delta$ be a weak solution to \eqref{eq:pde2o} with $m_\delta$ in place of $m$ and initial condition $u_0$ in the sense of Theorem \ref{thm:exist2obd}, for each $\delta\in (0,\bar\delta)$. Then, there exists a vanishing sequence $\delta_k\to 0$ and a map $u:\R_{\ge 0}\to\X$ such that for the sequence $(u_{\delta_k})_{k\in\N}$ and the limit $u$, one has
\begin{enumerate}[(a)]
\item $u\in L^\infty([0,T];L^2(\Omega))\cap L^2([0,T];H^1(\Omega))\cap C^{1/2}([0,T];\W_{\tilde m})$, and $u(t,x)\le S$ for almost all $t\ge 0$ and $x\in\Omega$;
\item $u_{\delta_k}(t)\to u(t)$ uniformly in $t\in[0,T]$ with respect to the distance $\W_{\tilde m}$ on $\X(0)$;
\item $u_{\delta_k}$ converges to $u$ strongly in $L^2([0,T];L^2(\Omega))$ and weakly in $L^2([0,T];H^1(\Omega))$;
\item $u$ is a solution to \eqref{eq:pde2o} in the sense of distributions, it attains the boundary conditions \eqref{eq:bc} and the initial condition $u(0,\cdot)=u_0$ holds a.e. on $\Omega$.
\end{enumerate}
\end{thm}
Again, the fourth-order case is similar:
\begin{thm}[Existence for non-Lipschitz mobilities: fourth order case (E2')]\label{thm:exist4onl}
Assume that the mobility $m$ satisfies \textup{(M1')--(M5')} and that the energy functional $\ent$ is of the form \textup{(E2')}. Let an initial datum $u_0\in \X\cap H^1(\Rd)$ with $u_0(x)\in [0,S]$ for almost every $x\in\Omega$ be given. Define $m_\delta$ for $\delta\in (0,\bar\delta)$ and sufficiently small $\bar\delta$ as in Definition \ref{def:appr} and let $u_\delta$ be a weak solution to \eqref{eq:pde4o} with $m_\delta$ in place of $m$ and initial condition $u_0$ in the sense of Theorem \ref{thm:exist4obd}, for each $\delta\in (0,\bar\delta)$. Then, there exists a vanishing sequence $\delta_k\to 0$ and a map $u:\R_{\ge 0}\to\X$ such that for the sequence $(u_{\delta_k})_{k\in\N}$ and the limit $u$, one has
\begin{enumerate}[(a)]
\item $u\in L^\infty([0,T];H^1(\Omega))\cap L^2([0,T];H^2(\Omega))\cap C^{1/2}([0,T];\W_{\tilde m})$, and $u(t,x)\le S$ for almost all $t\ge 0$ and $x\in\Omega$;
\item $u_{\delta_k}(t)\to u(t)$ uniformly in $t\in[0,T]$ with respect to the distance $\W_{\tilde m}$ on $\X(0)$;
\item $u_{\delta_k}$ converges to $u$ strongly in $L^2([0,T];H^1(\Omega))$ and weakly in $L^2([0,T];H^2(\Omega))$;
\item $u$ is a solution to \eqref{eq:pde4o} in the sense of distributions, it attains the boundary conditions \eqref{eq:bc} and the initial condition $u(0,\cdot)=u_0$ holds a.e. on $\Omega$.
\end{enumerate}
\end{thm}


\bibliographystyle{abbrv}
   \bibliography{bib}

\end{document}